\documentclass[12pt,twoside,a4paper]{article}
\usepackage[english]{babel}
\usepackage{graphicx} 
\usepackage{amssymb}
\usepackage[T1]{fontenc}
\usepackage{amsmath,epsf}
\newtheorem{Th}{Theorem}
\newtheorem{lemma}{Lemma}

\newcommand{\K}{\text{K}}
\newcommand{\I}{\text{I}}
\newcommand{\E}{\mathbb{E}}
\newcommand{\F}{\mathcal{F}}
\newcommand{\R}{\mathbb{R}}
\newcommand{\Z}{\mathbb{Z}}
\newcommand{\N}{\mathbb{N}}

\newcommand{\G}{\mathcal{G}}

\renewcommand{\P}{\mathbb{P}}

\newcounter{tictac}

\def\1{\,\rlap{\mbox{\small\rm 1}}\kern.15em 1}
\def\ind#1{\1_{#1}}
\def\build#1_#2^#3{\mathrel{\mathop{\kern 0pt#1}\limits_{#2}^{#3}}}
\def\tend#1#2{\build\hbox to 12mm{\rightarrowfill}_{#1\rightarrow #2}^{a.s.}}

\def\converge#1#2#3{\build\hbox to
15mm{\rightarrowfill}_{#1\rightarrow #2}^{\hbox{\scriptsize #3}}}
\begin{document}

\begin{center}
{\large
    {\sc
Kernel density estimation for stationary random fields
    }
}
\bigskip

Mohamed EL MACHKOURI

\medskip
\scriptsize{Laboratoire de Math\'ematiques Rapha\"el Salem\\
UMR CNRS 6085, Universit\'e de Rouen (France)\\
\emph{mohamed.elmachkouri@univ-rouen.fr}}
\end{center}

{\renewcommand\abstractname{Abstract}
\begin{abstract}
\baselineskip=18pt In this paper, under natural and easily verifiable conditions, we prove the $\mathbb{L}^1$-convergence and the 
asymptotic normality of the Parzen-Rosenblatt density estimator for stationary random fields of the form 
$X_k = g\left(\varepsilon_{k-s}, s \in \Z^d \right)$, $k\in\Z^d$, where $(\varepsilon_i)_{i\in\Z^d}$ are independent and identically distributed real 
random variables and $g$ is a measurable function defined on $\R^{\Z^d}$. Such kind of processes provides a general 
framework for stationary ergodic random fields. A Berry-Esseen's type central limit theorem is also given for 
the considered estimator.
{\em Short title:} Kernel density estimation for random fields.
\end{abstract}

\thispagestyle{empty}
\baselineskip=18pt
\section{Introduction and main results}
Let $(X_i)_{i\in\Z}$ be a stationary sequence of real random variables defined on a probability space $(\Omega, \F, \P)$ with an unknown marginal density $f$. 
The kernel density estimator $f_n$ of $f$ introduced by Rosenblatt \cite{Ros} and Parzen \cite{Parzen1962} is defined for all positive integer $n$ and any real $x$ by
$$
f_n(x)=\frac{1}{nb_n}\sum_{i=1}^n\K \left(\frac{x-X_i}{b_n}\right)
$$
where $\K$ is a probability kernel and the bandwidth $b_n$ is a parameter which converges slowly to zero such that $nb_n$ goes to infinity. 
The literature dealing with the asymptotic properties of $f_n$ when the observations $(X_i)_{i\in\Z}$ are independent is very extensive (see Silverman \cite{Silverman1986}). 
Parzen \cite{Parzen1962} proved that when $(X_i)_{i\in\Z}$ are independent and identically distribut (i.i.d) and the bandwidth $b_n$ goes to zero such that $nb_n$ goes to infinity 
then $(nb_n)^{1/2}(f_n(x_0)-\E f_n(x_0))$ converges in distribution to the normal law with zero mean and variance $f(x_0)\int_{\R} \K ^2(t)dt$. Under the 
same conditions on the bandwidth, this result was extended by Wu an Mielniczuk \cite{Wu-Mielniczuk} for causal linear processes with i.i.d. 
innovations and by Dedecker and Merlev\`ede \cite{Dedecker_Merlevede2002} for strongly mixing sequences.\\ 
In this paper, we are interested by the kernel 
density estimation problem in the setting of dependent random fields indexed by $\Z^d$ where $d$ is a positive integer. 
The question is not trivial since $\Z^d$ does not have a natural ordering for $d\geq 2$. In recent years, there is a growing interest in asymptotic properties 
of kernel density estimators for random fields. One can refer for example to Carbon et al. (\cite{Carbon-Hallin-Tran}, \cite{Carbon-Tran-Wu}), Cheng et al. 
\cite{Cheng-Ho-Lu2008}, El Machkouri \cite{Elmachkouri2011}, Hallin et al. \cite{Hallin--Lu--Tran2001}, Tran \cite{Tran} and Wang and Woodroofe \cite{Wang--Woodroofe2011}. 
In \cite{Tran}, the asymptotic normality of the kernel density estimator for strongly mixing random fields 
was obtained using the Bernstein's blocking technique and coupling arguments. Using the same method, the case of linear random fields with i.i.d. innovations 
was handled in \cite{Hallin--Lu--Tran2001}. In \cite{Elmachkouri2011}, the central limit theorem for the Parzen-Rosenblatt estimator given in \cite{Tran} 
was improved using the Lindeberg's method (see \cite{Lindeberg}) which seems to be better than the Bernstein's blocking technique approach. 
In particular, a simple criterion on the strong mixing coefficients is provided and the only 
condition imposed on the bandwith is $n^db_n\to\infty$ which is similar to the usual condition imposed in the independent case (see Parzen \cite{Parzen1962}). 
In \cite{Elmachkouri2011}, the regions where the random field is observed are reduced to squares but a carrefull reading of the proof allows us to 
state that the main result in \cite{Elmachkouri2011} still holds for very general regions $\Lambda_n$, namely those which the cardinality $\vert\Lambda_n\vert$ 
goes to infinity such that $\vert\Lambda_n\vert b_n$ goes to zero as $n$ goes to infinity (see Assumption $\textbf{(A3)}$ below). 
In \cite{Cheng-Ho-Lu2008}, Cheng et al. investigated the asymptotic normality of the kernel density estimator for linear random fields with i.i.d. innovations using a 
martingale approximation method (initiated by Cheng and Ho \cite{Cheng-Ho2006}) but it seems that there is a mistake in their proof (see Remark 6 in 
\cite{Wang--Woodroofe2011}). Since the mixing property is often unverifiable and might be too restrictive, it is important to provide limit theorems for 
nonmixing and possibly nonlinear random fields. We consider in this work a field $(X_i)_{i\in\Z^d}$ of identically distributed real random variables with an 
unknown marginal density $f$ such that 
\begin{equation}\label{definition_champ}
X_i = g\left(\varepsilon_{i-s};\,s\in\Z^d\right),\quad i\in\Z^d,
\end{equation}
where $(\varepsilon_j)_{j\in\Z^d}$ are i.i.d. random variables and $g$ is a measurable function defined on $\R^{\Z^d}$. 
In the one-dimensional case ($d=1$), the class (\ref{definition_champ}) includes linear as well as many widely used nonlinear time series models as special cases. More importantly, 
it provides a very general framework for asymptotic theory for statistics of stationary time series (see e.g. \cite{Wu2005} and the review paper \cite{Wu2011}).\\
We introduce the physical dependence measure first introduced by Wu \cite{Wu2005}. Let $(\varepsilon_j^{'})_{j\in\Z^d}$ be an i.i.d. copy of $(\varepsilon_j)_{j\in\Z^d}$ and consider for all positive integer $n$ the coupled version 
$X_i^{\ast}$ of $X_i$ defined by $X_i^{\ast}=g\left(\varepsilon^{\ast}_{i-s}\,;\,s\in\Z^d\right)$ where 
$\varepsilon_j^{\ast}=\varepsilon_j\ind{\{j\neq 0\}}+\varepsilon_0^{'}\ind{\{j=0\}}$ for all $j$ in $\Z^d$. In other words, we obtain $X_i^{\ast}$ from $X_i$ 
by just replacing $\varepsilon_0$ by its copy $\varepsilon_0^{'}$. Let $i$ in $\Z^d$ and $p>0$ be fixed. If $X_i$ belongs to $\mathbb{L}_{p}$ (that is, $\E\vert X_i\vert^p$ is finite), 
we define the physical dependence measure $\delta_{i,p}=\|X_i-X_i^{\ast}\|_{p}$ where $\|\,.\,\|_p$ is the usual $\mathbb{L}^p$-norm and we say that the random field 
$(X_i)_{i\in\Z^d}$ is $p$-stable if $\sum_{i\in\Z^d}\delta_{i,p}<\infty$. For $d\geq 2$, the reader should keep in mind the following two examples already given in \cite{Elmachkouri--Volny--Wu2013} :\\
\underline{{\em Linear random fields}}: Let $(\varepsilon_i)_{i\in\Z^d}$ be
i.i.d random variables with $\varepsilon_i$ in $\mathbb{L}^p$, $p
\geq 2$. The linear random field $X$ defined for all $i$ in $\Z^d$
by
$$
X_i=\sum_{s\in\Z^d}a_s\varepsilon_{i-s}
$$
with $(a_s)_{s\in\Z^d}$ in $\R^{\Z^d}$ such that $\sum_{i\in\Z^d}a_i^2<\infty$ is of the form $(\ref{definition_champ})$ with a linear functional
$g$. For all $i$ in $\Z^d$, $\delta_{i,p} = \vert a_i\vert \|\varepsilon_0 - \varepsilon^{'}_0\|_p$. So, $X$ is $p$-stable if $\sum_{i\in\Z^d}\vert a_i\vert<\infty$. 
Clearly, if $\textrm{H}$ is a Lipschitz continuous function, under the
above condition, the subordinated process $Y_i = \textrm{H}(X_i)$ is also $p$-stable since $\delta_{i,p} = O(|a_i|)$. \\
\underline{{\em Volterra field}} : Another class of nonlinear random field is the Volterra process which plays an important role in the 
nonlinear system theory (Casti \cite{Casti1985}, Rugh \cite{Rugh1981}): consider the second order Volterra process
\begin{eqnarray*}
X_i = \sum_{s_1, s_2\in\Z^d} a_{s_1, s_2} \varepsilon_{i-s_1} \varepsilon_{i-s_2},
\end{eqnarray*}
where $a_{s_1, s_2}$ are real coefficients with $a_{s_1, s_2} = 0$ if $s_1 = s_2$ and $(\varepsilon_i)_{i\in\Z^d}$ are i.i.d. random variables with $\varepsilon_i$ 
in $\mathbb{L}^p$, $p \geq 2$.
Let
\begin{eqnarray*}
A_i = \sum_{s_1, s_2\in\Z^d} (a_{s_1, i}^2 + a_{i, s_2}^2)\quad\textrm{and}\quad 
B_i = \sum_{s_1, s_2\in\Z^d} (|a_{s_1, i}|^p + |a_{i, s_2}|^p).
\end{eqnarray*}
By the Rosenthal inequality, there exists a constant $C_p > 0$ such that
\begin{eqnarray*}
\delta_{i,p} = \| X_i - X_i^*\|_p \leq C_p A_i^{1/2} \| \varepsilon_0\|_2 \| \varepsilon_0\|_p+ C_p B_i^{1/p} \| \varepsilon_0\|_p^2.
\end{eqnarray*}
From now on, for all finite subset $\Lambda$ of $\Z^d$, we denote $\vert\Lambda\vert$ the number of elements in $\Lambda$ and we 
observe $(X_i)_{i\in\Z^d}$ on a sequence $(\Lambda_n)_{n\geq 1}$ of finite subsets of $\Z^d$ which only satisfies $\vert\Lambda_n\vert$ 
goes to infinity as $n$ goes to infinity. It is important to note that we do not impose any condition on the boundary of the regions $\Lambda_n$. 
The density estimator $f_n$ of $f$ is defined for all positive integer $n$ and any real $x$ by
$$
f_n(x)=\frac{1}{\vert\Lambda_n\vert b_n}\sum_{i\in\Lambda_n}\K\left(\frac{x-X_i}{b_n}\right)
$$
where $b_n$ is the bandwidth parameter and $\K $ is a probability kernel. Our aim is to provide sufficient conditions for the $\mathbb{L}_1$-distance between 
$f_n$ and $f$ to converge to zero (Theorem \ref{L1_convergence}) and for 
$(\vert\Lambda_n\vert b_n)^{1/2}(f_n(x_i)-\E  f_n(x_i))_{1\leq i\leq k},\,(x_i)_{1\leq i\leq k}\in\R^k,\,k\in\N\backslash\{0\},$ to 
converge in law to a multivariate normal distribution (Theorem $\ref{convergence-loi}$) under minimal conditions on the bandwidth parameter. We give also a 
Berry-Esseen's type central limit theorem for the considered estimator (Theorem \ref{Berry-Esseen_type_clt}). 
In the sequel, we denote $\vert i\vert=\max_{1\leq k\leq d}\vert i_k\vert$ for all $i=(i_1,...,i_d)\in\Z^d$ and we denote also 
$\delta_i$ for $\delta_{i,2}$. The following assumptions are required.
\begin{itemize}
\item[\textbf{(A1)}] The marginal density function $f$ of each $X_k$ is Lipschitz.
\item[\textbf{(A2)}] $\K$ is Lipschitz, $\int_{\R}\K(u)\,du=1$, $\int_{\R}u^2\vert\K(u)\vert\,du<\infty$ and $\int_{\R}\K^2(u)\,du<\infty$.
\item[\textbf{(A3)}] $b_n\to0$ and $\vert\Lambda_n\vert\to\infty$ such that $\vert\Lambda_n\vert b_n\to\infty$.
\item[\textbf{(A4)}] $\sum_{i\in\Z^d}\vert i\vert^{\frac{5d}{2}}\,\delta_i<\infty$.
\end{itemize}
\begin{Th}\label{L1_convergence}
If $\textbf{\emph{(A1)}}$, $\textbf{\emph{(A2)}}$, $\textbf{\emph{(A3)}}$ and $\textbf{\emph{(A4)}}$ hold, then there exists $\kappa>0$ such that for all 
integer $n\geq 1$,
\begin{equation}\label{convergence_L1_with_rate}
\E\int_\R\vert f_n(x)-f(x)\vert\,dx\leq\kappa\left(b_n+\frac{1}{\sqrt{\vert\Lambda_n\vert b_n}}\right)^{\frac{2}{3}}.
\end{equation}
\end{Th}
\textbf{Remark 1.} One can optimize the inequality (\ref{convergence_L1_with_rate}) by taking $b_n=\vert\Lambda_n\vert^{-\frac{1}{3}}$. Then, we obtain 
$\E\int_\R\vert f_n(x)-f(x)\vert\,dx=O\left(\vert\Lambda_n\vert^{-\frac{2}{9}}\right)$.\\
\\
\textbf{Remark 2.} The convergence in probability of $\int_\R\vert f_n(x)-f(x)\vert\,dx$ to $0$ was obtained (without rate) by Hallin et al. 
(\cite{Hallin-Lu-Tran-2004a}, Theorem 2.1) for rectangular region $\Lambda_n$. The authors defined the so-called stability coefficients 
$(v(m))_{m\geq 1}$ by $v(m)=\|X_0-\overline{X}_0\|_2^2$ where $\overline{X}_0=\E\left(X_0\vert\mathcal{H}_m\right)$ and $\mathcal{H}_m=\sigma\left(\varepsilon_{s}\,,\,\vert s\vert\leq m\right)$. 
Under minimal conditions on the bandwidth $b_n$, with our notations, their result holds as soon as $v(m)=o(m^{-4d})$. Arguing as in the proof of Lemma 
\ref{controle_norme_de_la_difference_de_K_0_et_K_0_barre} below, one can relate the stability coefficients with the physical dependence measure ones by the inequality 
$v(m)\leq C\sum_{\vert i\vert>m}\delta_i^2$, $m\geq 1$, $C>0$.\\
\\
In the sequel, we consider the sequence $(m_n)_{n\geq 1}$ defined by 
\begin{equation}\label{def_mn}
m_n=\max\left\{v_n,\left[\left(\frac{1}{b_n^{3}}\sum_{\vert i\vert>v_n}\vert i\vert^{\frac{5d}{2}}\,\delta_i\right)^{\frac{1}{3d}}\right]+1\right\}
\end{equation}
where $v_n=\big[b_n^{-\frac{1}{2d}}\big]$ and $[\,.\,]$ denotes the integer part function. The following technical lemma is a spatial version of a 
result by Bosq et al. (\cite{Bosq--Merlevede--Peligrad1999}, pages 88-89). 
\begin{lemma}\label{mn}
If $\emph{\textbf{(A4)}}$ holds then   
$$
m_n\to\infty,\quad m_n^db_n\to0\quad\textrm{and}\quad\frac{1}{(m_n^db_n)^{3/2}}\sum_{\vert i\vert>m_n}\vert i\vert^{\frac{5d}{2}}\,\delta_i\to 0.
$$
\end{lemma}
For all $z$ in $\R$ and all $i$ in $\Z^d$, we denote
\begin{equation}\label{definition_K_i_z_overline_K_i_z}
\K_i(z)=\K\left(\frac{z-X_i}{b_n}\right)\quad\textrm{and}\quad\overline{\K}_i(z)=\E\left(\K_i(z)\vert\F_{n,i}\right)
\end{equation}
where $\F_{n,i}=\sigma\left(\varepsilon_{i-s}\,;\,\vert s\vert\leq m_n\right)$. So, denoting $M_n=2m_n+1$, $(\overline{\K}_i(z))_{i\in\Z^d}$ is 
an $M_n$-dependent random field (i.e. $\overline{\K}_i(z)$ and $\overline{\K}_j(z)$ are independent as soon as $\vert i-j\vert\geq M_n$).
\begin{lemma}\label{moment_inequality}
For all $p>1$, all $x$ in $\R$, all positive integer $n$ and all $(a_i)_{i\in\Z^d}$ in $\R^{\Z^d}$, 
$$
\left\|\sum_{i\in\Lambda_n}a_i\left(\emph{\K}_i(x)-\overline{\emph{\K}}_i(x)\right)\right\|_p
\leq\frac{8 m_n^d}{b_n}\left(p\sum_{i\in\Lambda_n}a_i^2\right)^{1/2}\sum_{\vert i\vert>m_n}\delta_{i,p}.
$$
\end{lemma} 
In order to establish the asymptotic normality of $f_n$, we need additional assumptions:
\begin{itemize}
\item[\textbf{(B1)}] The marginal density function of each $X_k$ is positive, continuous and bounded.
\item[\textbf{(B2)}] $\K$ is Lipschitz, $\int_{\R}\K(u)\,du=1$, $\int_{\R}\vert\K(u)\vert\,du<\infty$ and $\int_{\R}\K^2(u)\,du<\infty$.
\item[\textbf{(B3)}] There exists $\kappa>0$ such that $\sup_{\substack{(x,y)\in\R^2\\ i\in\Z^d\backslash\{0\}}} f_{0,i}(x,y)\leq\kappa$ where $f_{0,i}$ is the joint density of $(X_0,X_i)$.
\end{itemize}
\begin{Th}\label{convergence-loi}
Assume that $\textbf{\emph{(A3)}}$, $\textbf{\emph{(A4)}}$, $\textbf{\emph{(B1)}}$, $\textbf{\emph{(B2)}}$ and $\textbf{\emph{(B3)}}$ hold. 
For all positive integer $k$ and any distinct points $x_1,...,x_k$ in $\R$,
\begin{equation}\label{limit}
(\vert\Lambda_n\vert b_{n})^{1/2}
\left(\begin{array}{c}
       f_{n}(x_1)-\E f_n(x_1)\\
       \vdots\\
       f_{n}(x_k)-\E f_n(x_k)
       \end{array} \right)
\converge{n}{\infty}{\emph{Law}}
\mathcal{N}\left(0,\Gamma\right)
\end{equation}
where $\Gamma$ is a diagonal matrix with diagonal elements $\gamma_{ii}=f(x_i)\int_{\R}\emph{\K}^2(u)du$.
\end{Th}
\textbf{Remark 3.} A replacement of $\E f_n(x_i)$ by $f(x_i)$ for all $1\leq i\leq k$ in ($\ref{limit}$) 
is a classical problem in density estimation theory. Let $s\geq 2$ be a positive integer and $\kappa>0$. If the $s$th derivative $f^{(s)}$ of $f$ exists such that $\vert f^{(s)}\vert\leq\kappa$ and the kernel $\K$ satisfies $\int_{\R} u^r \K(u) du=0$ for $r=1,2,...,s-1$ and $0<\int_{\R} \vert u\vert^s \vert \K(u)\vert du<\infty$ then $\vert \E f_n(x_i)-f(x_i)\vert=O(b_n^s)$ and thus the 
centering $\E f_n(x_i)$ may be changed to $f(x_i)$ without affecting the above result provided that $\vert\Lambda_n\vert b_n^{2s+1}$ converges to zero.\\
\\
\textbf{Remark 4.} If $(X_i)_{i\in\Z^d}$ is a linear random field of the form $X_i=\sum_{j\in\Z^d}a_j\varepsilon_{i-j}$ where $(a_j)_{j\in\Z^d}$ 
are real numbers such that $\sum_{j\in\Z^d}a_j^2<\infty$ and $(\varepsilon_j)_{j\in\Z^d}$ are i.i.d. real random variables with zero mean and finite variance then 
$\delta_i=\vert a_i\vert\|\varepsilon_0-\varepsilon_0^{'}\|_2$ and Theorem \ref{convergence-loi} holds provided that 
$\sum_{i\in\Z^d}\vert i\vert^{\frac{5d}{2}}\vert a_i\vert<\infty$. For $\Lambda_n$ rectangular, Hallin et al. \cite{Hallin--Lu--Tran2001} obtained the same result when 
$\vert a_j\vert =O\left(\vert j\vert^{-\gamma}\right)$ with $\gamma>\max\{d+3,2d+0.5\}$ and $\vert\Lambda_n\vert b_n^{(2\gamma-1+6d)/(2\gamma-1-4d)}$ goes to infinity. So, in the particular 
case of linear random fields, our assumption $\textbf{(A4)}$ is more restrictive than the condition obtained by Hallin et al. \cite{Hallin--Lu--Tran2001} but 
our result is valid for a larger class of random fields and under only minimal conditions on the bandwidth (see Assumption \textbf{(A3)}). 
Finally, for causal linear random fields, 
Wang and Woodroofe \cite{Wang--Woodroofe2011} obtained also a sufficient condition on the coefficients $(a_j)_{j\in\N^d}$ for the kernel density estimator to be 
asymptotically normal. Their condition is less restrictive than the condition $\sum_{i\in\Z^d}\vert i\vert^{\frac{5d}{2}}\vert a_i\vert<\infty$ but they assumed also 
$\E(\vert\varepsilon_0\vert^{p})<\infty$ for some $p>2$.\\
\\
Now, we are going to investigate the rate of convergence in (\ref{limit}). For all positive integer $n$ and all $x$ in $\R$, we denote 
$\textrm{D}_n(x)=\sup_{t\in\R}\left\vert\P\left(U_n(x)\leq t\right)-\Phi(t)\right\vert$ where $\Phi$ is the distribution function of the standard normal law and 
$$
U_n(x)=\frac{\sqrt{\vert\Lambda_n\vert b_n}\left(f_n(x)-\E f_n(x)\right)}{\sqrt{f(x)\int_{\R}\K^2(t)dt}}.
$$
\begin{Th}\label{Berry-Esseen_type_clt}
Let $n$ in $\N\backslash\{0\}$ and $x$ in $\R$ be fixed. Assume that $\int_{\R}\vert \emph{K}(t)\vert^\tau dt<\infty$ for some $2<\tau\leq 3$. If there exist $\alpha>1$ and $p\geq 2$ such that $\sum_{i\in\Z^d}\vert i\vert^{d\alpha}\delta_{i,p}<\infty$ then there exists a constant $\kappa>0$ such that $\emph{D}_n(x)\leq\kappa\vert\Lambda_n\vert^{-\theta}$ where
$$
\theta=\theta(\alpha,\tau,p)=\left(\frac{1}{2}-\frac{1}{\tau}\right)\frac{3p(1-\tau)+2p(\alpha-1)}{(\tau-1)(p+1)+p(\alpha-1)}.
$$
\end{Th}
\textbf{Remark 5}. If $\tau=3$, $p=2$ and $\sum_{i\in\Z^d}\vert i\vert^{d\alpha}\delta_{i}<\infty$ for some $\alpha>4$ then 
$$
\emph{D}_n(x)\leq\kappa\vert\Lambda_n\vert^{-\theta(\alpha)}
\quad\textrm{where}\quad
\theta(\alpha)=\frac{2\alpha-8}{3(4+2\alpha)}\converge{\alpha}{\infty}{}\frac{1}{3}.
$$
\section{Numerical illustration}
In this section, we give some simulations with a view to illustrate the results given in this paper. We assume $d=2$ and we consider the autoregressive random field $(X_{i,j})_{(i,j)\in\Z^2}$ defined by 
\begin{equation}\label{model}
X_{i,j}=\alpha X_{i-1,j}+\beta X_{i,j-1}+\varepsilon_{i,j}
\end{equation}
where $\alpha=0.2$, $\beta=0.7$ and $(\varepsilon_{i,j})_{(i,j)\in\Z^2}$ are iid random variables uniformly distributed over the intervalle $[-5,5]$. Since $\vert \alpha\vert+\vert \beta\vert<1$, the equation $(\ref{model})$ has a stationary solution $X_{i,j}$ (see \cite{Kulkarni1992}) defined by
\begin{equation}\label{model_bis}
X_{i,j}=\sum_{k_1\geq 0}\sum_{k_2\geq 0} {k_1+k_2 \choose k_1}\alpha^{k_1}\beta^{k_2}\varepsilon_{i-k_1,j-k_2}
\end{equation}
and each $X_{i,j}$ is uniformly distributed over the intervalle $[-5\gamma,5\gamma]$ with 
$$
\gamma=\sum_{k_1\geq 0}\sum_{k_2\geq 0} {k_1+k_2 \choose k_1}\alpha^{k_1}\beta^{k_2}=\frac{1}{1-(\alpha+\beta)}=10.
$$
We simulate the $\varepsilon_{i,j}$'s over the rectangular grid $[0,2t]^2\cap\Z^2$ where $t$ is a positive integer and the data $X_{i,j}$ over the grid $\Lambda_t=[t+1,2t]^2\cap\Z^2$ following (\ref{model_bis}). We take the data $X_{i,j}$ for $(i,j)$ in the region $\Lambda_t$ as our data set and we calculate from this data set the kernel density estimator
\begin{equation}\label{kernel_density_estimator}
\hat{f}_t(x)=\frac{1}{t^2\times b_t}\sum_{(i,j)\in\Lambda_t}\textrm{K}\left(\frac{x-X_{i,j}}{b_t}\right)
\end{equation}
where $x$ is fixed in $\R$, $b_t$ is the bandwith parameter and $\textrm{K}$ is the Epanachnikov kernel defined by $\textrm{K}(s)=\frac{3}{4}(1-s^2)$ if $s\in]-1,1[$ and $\textrm{K}(s)=0$ if $s\notin]-1,1[$.\\
In order to illustrate the result obtained in Theorem \ref{L1_convergence}, we calculate (Monte Carlo method) $\int_{-100}^{100}\vert\hat{f}_t(x)-f(x)\vert dx$ where $f$ is the true density function of $X_{0,0}$ and the bandwith $b_t$ is being set to $\vert\Lambda_t\vert^{-1/3}$ with $\vert\Lambda_t\vert$ denoting the number of elements in $\Lambda_t$. Hence, we derive its expectation $\E\int_{-100}^{100}\vert\hat{f}_t(x)-f(x)\vert dx$ by taking the arithmetic mean value of $100$ replications of $\int_{-100}^{100}\vert\hat{f}_t(x)-f(x)\vert dx$. The results are given for several values of $t$ in the following table
\begin{center}
\begin{tabular}{|c|c|c|c|}
\hline
$t$ &  $\vert\Lambda_t\vert=t^2$ & $b_t=\vert\Lambda_t\vert^{-1/3}$ & $\E\int_{-100}^{100}\vert\hat{f}_t(x)-f(x)\vert dx$\\
\hline\hline
$10$ & $100$ & $0.215$ & $0.0171$\\
\hline
$20$ & $400$ & $0.136$ & $0.0163$\\
\hline
$50$ & $2500$ & $0.074$ & $0.0157$\\
\hline
$100$ & $10000$ & $0.046$ & $0.0153$\\
\hline
\end{tabular}
\end{center}
and we observe the $L^1$-convergence of $\hat{f}_t$ to the true density function $f$ of $X_{0,0}$. In order to illustrate the asymptotic normality of the estimator (\ref{kernel_density_estimator}), we put $x=-1$, $t=20$ and $b_{20}=0.7$ and we calculate the expectation $\E\left(\hat{f}_t(-1)\right)$ of $\hat{f}_t(-1)$ by taking again the arithmetic mean value of $100$ replications of $\hat{f}_t(-1)$. Finally, noting that $\int_{\R}\textrm{K}^2(x)dx=4/5$ and $f(-1)=1/100$, we consider $1500$ replications of 
$$
\frac{\sqrt{400\times 0.7}\left(\hat{f}_{20}(-1)-\E\left(\hat{f}_{20}(-1)\right)\right)}{\sqrt{1/100\times 4/5}}
$$
and we obtain the following histogram (see figure \ref{Histogram}) which seems to fit well to the target distribution, that is the standard normal law $\mathcal{N}(0,1)$.
\begin{figure}[!h] 
\begin{center}
\includegraphics[width=10cm]{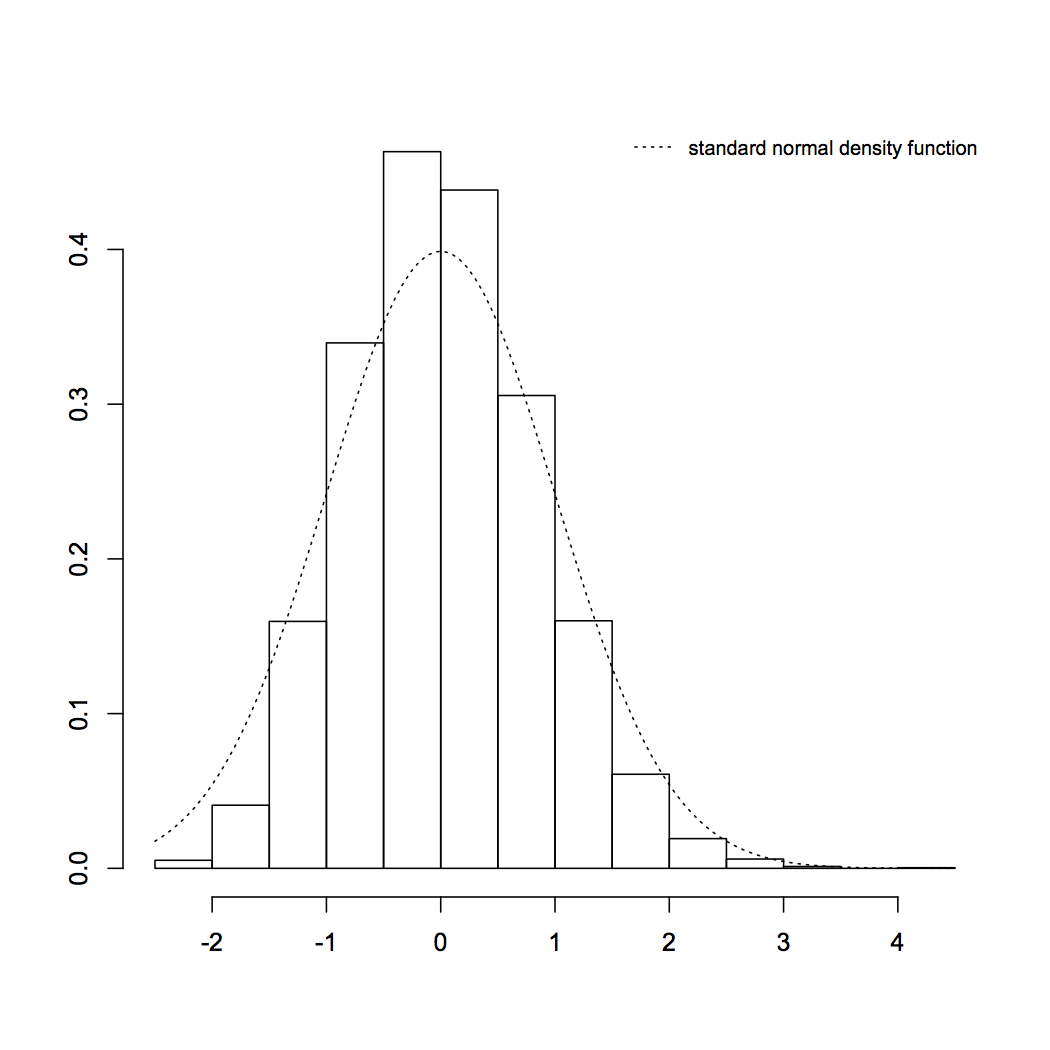} 
\end{center}
\caption{Asymptotic normality of the kernel density estimator.} 
\label{Histogram} 
\end{figure} 

In the simulation given in Figure \ref{Histogram}, we fixed the bandwith $b_{20}=0.7$ arbitrarily since we do not investigate in this work any procedure for a data-driven choice of the bandwith parameter. Such a study is an important task and will be done in a forthcoming paper.

\section{Proofs}
The proof of all lemmas of this section are postponed to the appendix. In the sequel, the letter $\kappa$ denotes a positive constant which the value is not important. 
\subsection{Proof of Theorem $\textbf{\ref{L1_convergence}}$}
For all positive integer $n$, denote $\textrm{J}_n=\int_\R\vert f_n(x)-f(x)\vert\,dx$. For all real $A\geq 1$, we have 
$\textrm{J}_n=\textrm{J}_{n,1}(A)+\textrm{J}_{n,2}(A)$ where 
$$
\textrm{J}_{n,1}(A)=\int_{\vert x\vert>A}\vert f_n(x)-f(x)\vert\,dx\quad\textrm{and}\quad \textrm{J}_{n,2}(A)=\int_{\vert x\vert\leq A}\vert f_n(x)-f(x)\vert\,dx.
$$
Moreover
$$
\E\textrm{J}_{n,1}(A)\leq\int_{\vert x\vert>A}\E\vert f_n(x)\vert dx+\frac{1}{A^2}\int_{\R}x^2 f(x)dx
$$
and
\begin{align*}
\int_{\vert x\vert>A}\E\vert f_n(x)\vert dx&\leq\int_{\vert x\vert>A}\int_{\R}\vert\textrm{K}(t)\vert f(x-b_nt)dtdx\\
&=\int_{\vert t\vert>\frac{A}{2}}\vert\textrm{K}(t)\vert\int_{\vert x\vert>A}f(x-b_nt)dxdt+\int_{\vert t\vert\leq\frac{A}{2}}\vert\textrm{K}(t)\vert\int_{\vert x\vert>A}f(x-b_nt)dxdt\\
&\leq \int_{\vert t\vert>\frac{A}{2}}\vert\textrm{K}(t)\vert\int_{\vert y+b_nt\vert>A}f(y)dydt+\int_{\vert t\vert\leq\frac{A}{2}}\vert\textrm{K}(t)\vert\int_{\vert y\vert>A(1-\frac{b_n}{2})}f(y)dydt\\
&\leq\frac{4}{A^2}\int_{\R}t^2\vert\textrm{K}(t)\vert dt+\frac{4}{A^2}\int_{\R}\vert\textrm{K}(t)\vert dt\int_{\R}y^2 f(y)dy. 
\end{align*}
Consequently, we obtain 
\begin{equation}\label{majoration_EJn1A}
\E\textrm{J}_{n,1}(A)\leq\frac{\kappa}{A^2}.
\end{equation}
Now, $\textrm{J}_{n,2}(A)\leq \textrm{J}_{n,2}^{(1)}(A)+\textrm{J}_{n,2}^{(2)}(A)$ where
$$
\textrm{J}_{n,2}^{(1)}(A)=\int_{\vert x\vert\leq A}\vert f_n(x)-\E f_n(x)\vert\,dx\quad\textrm{and}\quad 
\textrm{J}_{n,2}^{(2)}(A)=\int_{\vert x\vert\leq A}\vert\E f_n(x)-f(x)\vert\,dx.
$$
Since 
\begin{align*}
\vert\E f_n(x)-f(x)\vert&=\left\vert\int_{\R}\textrm{K}(t)\left(f(x-b_n t)-f(x)\right)dt\right\vert\\
&\leq\int_{\R}\vert\textrm{K}(t)\vert\left\vert f(x-b_n t)-f(x)\right\vert dt\\
&\leq \kappa b_n\int_{\R}\vert t\vert \vert\textrm{K}(t)\vert dt,
\end{align*}
we obtain
\begin{equation}\label{majoration_Jn2A}
\textrm{J}_{n,2}^{(2)}(A)\leq \kappa A b_n.
\end{equation}
Keeping in mind the notation (\ref{definition_K_i_z_overline_K_i_z}) and denoting 
$\overline{f}_n(x)=\frac{1}{\vert\Lambda_n\vert b_n}\sum_{i\in\Lambda_n}\overline{\K}_i(x)$, we have $\textrm{J}_{n,2}^{(1)}(A)\leq\I_{n,1}(A)+\I_{n,2}(A)$ where 
$$
\I_{n,1}(A)=\int_{\vert x\vert\leq A}\vert f_n(x)-\overline{f}_n(x)\vert\,dx\quad\textrm{and}\quad
\I_{n,2}(A)=\int_{\vert x\vert\leq A}\vert \overline{f}_n(x)-\E\overline{f}_n(x)\vert\,dx.
$$
By Lemma \ref{moment_inequality}, we have
$$
\left\|f_n(x)-\overline{f}_n(x)\right\|_2\leq\frac{\kappa\sum_{\vert i\vert>m_n}\vert i\vert^{\frac{5d}{2}}\delta_i}{\sqrt{\vert\Lambda_n\vert b_n}(m_n^db_n)^{3/2}}.
$$
Applying Lemma \ref{mn}, we obtain 
\begin{equation}\label{majoration_EIn1A}
\E\I_{n,1}(A)\leq\frac{\kappa A}{\sqrt{\vert\Lambda_n\vert b_n}}. 
\end{equation}
Now, $\left\|\overline{f}_n(x)-\E\overline{f}_n(x)\right\|_2^2$ equals to
\begin{equation}\label{developpement_norme_2_au_carre}
\frac{1}{\vert\Lambda_n\vert^{2}b_n}\left(\vert\Lambda_n\vert\E\left(\overline{Z}_0^2(x)\right)
+\sum_{\substack{j\in\Z^d\backslash\{0\} \\ \vert j\vert< M_n}}\vert\Lambda_n\cap(\Lambda_n-j)\vert\E\left(\overline{Z}_0(x)\overline{Z}_j(x)\right)\right)
\end{equation}
where we recall that $\overline{Z}_i(x)=\frac{1}{\sqrt{b_n}}\left(\overline{\K}_i(x)-\E\overline{\K}_i(x)\right)$ and $M_n=2m_n+1$.
\begin{lemma}\label{lemme_technique_0}
Let $x$, $s$ and $t$ be fixed in $\R$. Then $\E\left(\overline{Z}_0^2(x)\right)$ converges to $f(x)\int_{\R}\emph{\K}^2(u)du$ and 
$\sup_{i\in\Z^d\backslash\{0\}}\E\vert\overline{Z}_0(s)\overline{Z}_i(t)\vert=o(M_n^{-d})$.
\end{lemma}
Combining (\ref{developpement_norme_2_au_carre}) and Lemma \ref{lemme_technique_0}, we derive 
$\left\|\overline{f}_n(x)-\E\overline{f}_n(x)\right\|_2^2=O\left(\left(\vert\Lambda_n\vert b_n\right)^{-1}\right)$. Hence, 
\begin{equation}\label{majoration_EIn2A}
\E\I_{n,2}(A)\leq\frac{\kappa A}{\sqrt{\vert\Lambda_n\vert b_n}}. 
\end{equation}
Combining (\ref{majoration_EJn1A}), (\ref{majoration_Jn2A}), (\ref{majoration_EIn1A}) and (\ref{majoration_EIn2A}), we obtain 
$$
\E\textrm{J}_n\leq\kappa\left(\frac{1}{A^2}+A\left(b_n+\frac{1}{\sqrt{\vert\Lambda_n\vert b_n}}\right)\right).
$$
Optimizing in $A$, we derive ($\ref{convergence_L1_with_rate}$). The proof of Theorem \ref{L1_convergence} is complete.
\subsection{Proof of Theorem $\textbf{\ref{convergence-loi}}$} 
Without loss of generality, we consider only the case $k=2$ and we refer to $x_1$ and $x_2$ as $x$ and $y$ ($x\neq y$). 
Let $\lambda_1$ and $\lambda_2$ be two constants such that $\lambda_1^2+\lambda_2^2=1$ and note that 
\begin{align*}
\lambda_1(\vert\Lambda_n\vert b_n)^{1/2}(f_n(x)-\E f_n(x))+\lambda_2(\vert\Lambda_n\vert b_n)^{1/2}(f_n(y)-\E f_n(y))
&=\sum_{i\in\Lambda_n}\frac{\Delta_i}{\vert\Lambda_n\vert^{1/2}},\\
\lambda_1(\vert\Lambda_n\vert b_n)^{1/2}(\overline{f}_n(x)-\E\overline{f}_n(x))+\lambda_2(\vert\Lambda_n\vert b_n)^{1/2}(\overline{f}_n(y)-\E\overline{f}_n(y))
&=\sum_{i\in\Lambda_n}\frac{\overline{\Delta}_i}{\vert\Lambda_n\vert^{1/2}},
\end{align*} 
where $\Delta_i=\lambda_1 Z_{i}(x)+\lambda_2 Z_{i}(y)$ and $\overline{\Delta}_i=\lambda_1\overline{Z}_{i}(x)+\lambda_2\overline{Z}_{i}(y)$ and for all $z$ in $\R$,
$$
Z_{i}(z)=\frac{1}{\sqrt{b_n}}\left(\K_i(z)-\E\K_i(z)\right)
\quad\textrm{and}\quad
\overline{Z}_{i}(z)=\frac{1}{\sqrt{b_n}}\left(\overline{\K}_i(z)-\E\overline{\K}_i(z)\right)
$$
where $\K_i(z)$ and $\overline{\K}_i(z)$ are defined by $(\ref{definition_K_i_z_overline_K_i_z})$. Applying Lemma \ref{mn} and Lemma \ref{moment_inequality}, we know that
\begin{equation}\label{ecart_sommes_partielles_delta_et_overline_delta}
\frac{1}{\vert\Lambda_n\vert^{1/2}}\left\|\sum_{i\in\Lambda_n}\left(\Delta_i-\overline{\Delta}_i\right)\right\|_2
\leq\frac{\kappa(\vert\lambda_1\vert+\vert\lambda_2\vert)}{(m_n^db_n)^{3/2}}\sum_{\vert i\vert>m_n}\vert i\vert^{\frac{5d}{2}}\delta_i=o(1).
\end{equation}
So, it suffices to prove the asymptotic normality of the sequence $\left(\vert\Lambda_n\vert^{-1/2}\sum_{i\in\Lambda_n}\overline{\Delta}_i\right)_{n\geq 1}$. 
We are going to follow the Lindeberg's type proof of Theorem $1$ in \cite{Dedecker1998}. We consider the notations
\begin{equation}\label{eta_sigma}
\eta=(\lambda_1^2f(x)+\lambda_2^2f(y))\sigma^2\quad\textrm{and}\quad\sigma^2=\int_{\R} \K ^2(u)du.
\end{equation}
\begin{lemma}\label{lemme-technique} $\E(\overline{\Delta}_0^2)$ converges to $\eta$ and 
$\sup_{i\in\Z^d\backslash\{0\}}\E\vert\overline{\Delta}_0\overline{\Delta}_i\vert=o(M_n^{-d})$.
\end{lemma}
On the lattice $\Z^{d}$ we define the lexicographic order as follows: if
$i=(i_{1},...,i_{d})$ and $j=(j_{1},...,j_{d})$ are distinct elements of $\Z^{d}$, the notation $i<_{\textrm{lex}}j$ means that either
$i_{1}<j_{1}$ or for some $k$ in $\{2,3,...,d\}$, $i_{k}<j_{k}$ and $i_{l}=j_{l}$ for $1\leq l<k$. We let $\varphi$ denote the unique function from 
$\{1,...,\vert\Lambda_n\vert\}$ to $\Lambda_{n}$ such that $\varphi(k)<_{\text{\text{lex}}}\varphi(l)$ 
for $1\leq k<l\leq\vert\Lambda_n\vert$. For all real random field $(\zeta_i)_{i\in\Z^d}$ and all integer $k$ in $\{1,...,\vert\Lambda_n\vert\}$, we denote
$$
S_{\varphi(k)}(\zeta)=\sum_{i=1}^k \zeta_{\varphi(i)}\quad\textrm{and}\quad
S_{\varphi(k)}^{c}(\zeta)=\sum_{i=k}^{\vert\Lambda_n\vert} \zeta_{\varphi(i)}
$$
with the convention $S_{\varphi(0)}(\zeta)=S_{\varphi(\vert\Lambda_n\vert+1)}^{c}(\zeta)=0$. From now on, we consider a field $(\xi_{i})_{i\in\Z^d}$ of i.i.d. 
standard normal random variables independent of $(X_{i})_{i\in\Z^d}$. We introduce the fields $Y$ and $\gamma$ defined for all $i$ in $\Z^d$ by 
$$
Y_{i}=\frac{\overline{\Delta}_i}{\vert\Lambda_n\vert^{1/2}}\quad\textrm{and}\quad\gamma_{i}=\frac{\sqrt{\eta}\xi_i}{\vert\Lambda_n\vert^{1/2}}
$$ 
where $\eta$ is defined by ($\ref{eta_sigma}$). Note that $Y$ is an $M_n$-dependent random field where $M_n=2m_n+1$ and $m_n$ is defined by (\ref{def_mn}). 
Let $h$ be any function from $\R$ to $\R$. For $0<k\leq l\leq \vert\Lambda_n\vert$, we introduce
$h_{k,l}(Y)=h(S_{\varphi(k)}(Y)+S_{\varphi(l)}^{c}(\gamma))$. With the above
convention we have that $h_{k,\vert\Lambda_n\vert+1}(Y)=h(S_{\varphi(k)}(Y))$ and also
$h_{0,l}(Y)=h(S_{\varphi(l)}^{c}(\gamma))$. In the sequel, we will often
write $h_{k,l}$ instead of $h_{k,l}(Y)$. We denote by $B_{1}^4(\R)$ the unit ball of $C_{b}^4(\R)$: $h$ belongs to
$B_{1}^4(\R)$ if and only if it belongs to $C^4(\R)$ and satisfies $\max_{0\leq i\leq 4}\|h^{(i)}\|_{\infty}\leq 1$. It suffices to prove that for all $h$ in $B_{1}^4(\R)$,
$$
\E\left(h\left(S_{\varphi(\vert\Lambda_n\vert)}(Y)\right)\right)\converge{n}{\infty}{}\E \left(h\left(\sqrt{\eta}\xi_0\right)\right).
$$
We use Lindeberg's decomposition:
$$
\E \left(h\left(S_{\varphi(\vert\Lambda_n\vert)}(Y)\right)-h\left(\sqrt{\eta}\xi_{0}\right)\right)
=\sum_{k=1}^{\vert\Lambda_n\vert}\E \left(h_{k,k+1}-h_{k-1,k}\right).
$$
Now, we have $h_{k,k+1}-h_{k-1,k}=h_{k,k+1}-h_{k-1,k+1}+h_{k-1,k+1}-h_{k-1,k}$ and by Taylor's formula we obtain
\begin{align*}
h_{k,k+1}-h_{k-1,k+1}&=Y_{\varphi(k)}h_{k-1,k+1}^{'}+\frac{1}{2}Y_{\varphi(k)}^{2}h_{k-1,k+1}^{''}+R_{k}\\
h_{k-1,k+1}-h_{k-1,k}&=-\gamma_{\varphi(k)}h_{k-1,k+1}^{'}-\frac{1}{2}\gamma_{\varphi(k)}^{2}h_{k-1,k+1}^{''}+r_{k}
\end{align*}
where $\vert R_{k}\vert\leq Y_{\varphi(k)}^2(1\wedge\vert Y_{\varphi(k)}\vert)$ and $\vert r_{k}\vert\leq\gamma_{\varphi(k)}^2(1\wedge\vert\gamma_{\varphi(k)}\vert)$.
Since $(Y,\xi_{i})_{i\neq \varphi(k)}$ is independent of $\xi_{\varphi(k)}$, it follows that
$$
\E \left(\gamma_{\varphi(k)}h_{k-1,k+1}^{'}\right)=0\quad\textrm{and}\quad
\E \left(\gamma_{\varphi(k)}^2h_{k-1,k+1}^{''}\right)=\E \left(\frac{\eta}{\vert\Lambda_n\vert}h_{k-1,k+1}^{''}\right)
$$
Hence, we obtain
\begin{align*}
\E \left(h(S_{\varphi(\vert\Lambda_n\vert)}(Y))-h\left(\sqrt{\eta}\xi_0\right)\right)&=
\sum_{k=1}^{\vert\Lambda_n\vert}\E (Y_{\varphi(k)}h_{k-1,k+1}^{'})\\
&\quad+\sum_{k=1}^{\vert\Lambda_n\vert}\E \left(\left(Y_{\varphi(k)}^2-\frac{\eta}{\vert\Lambda_n\vert}\right)\frac{h_{k-1,k+1}^{''}}{2}\right)\\
&\quad+\sum_{k=1}^{\vert\Lambda_n\vert}\E \left(R_{k}+r_{k}\right).
\end{align*}
Let $1\leq k\leq \vert\Lambda_n\vert$ be fixed. Since $\E\vert\overline{\Delta}_0\vert=O\left(\sqrt{b_n}\right)$ and 
$\left(\overline{\Delta}_0^2b_n\right)_{n\geq 1}$ is uniformly integrable, we derive
$$
\sum_{k=1}^{\vert\Lambda_n\vert}\E\vert R_k\vert\leq\E\left(\overline{\Delta}_0^2\left(1\wedge\frac{\vert\overline{\Delta}_0\vert}{\vert\Lambda_n\vert^{1/2}}\right)\right)
=o(1)
$$
and
$$
\sum_{k=1}^{\vert\Lambda_n\vert}\E\vert r_k\vert\leq\frac{\eta^{3/2}\E\vert\xi_0\vert^3}{\vert\Lambda_n\vert^{1/2}}
=O\left(\vert\Lambda_n\vert^{-1/2}\right).
$$
Consequently, we obtain
$$
\sum_{k=1}^{\vert\Lambda_n\vert}\E \left(\vert R_{k}\vert+\vert r_{k}\vert\right)=o(1).
$$
Now, it is sufficient to show
\begin{equation}\label{equation1}
\lim_{n\to\infty}\sum_{k=1}^{\vert\Lambda_n\vert}\left(\E (Y_{\varphi(k)}h_{k-1,k+1}^{'})+\E \left(\left(Y_{\varphi(k)}^2-\frac{\eta}{\vert\Lambda_n\vert}\right)\frac{h_{k-1,k+1}^{''}}{2}\right)\right)=0.
\end{equation}
First, we focus on $\sum_{k=1}^{\vert\Lambda_n\vert}\E \left(Y_{\varphi(k)}h_{k-1,k+1}^{'}\right)$. 
Let the sets $\{V_{i}^{k}\,;\,i\in\Z^{d}\,,\,k\in\N\backslash\{0\}\}$ be defined as follows: $V_{i}^{1}=\{j\in\Z^{d}\,;\,j<_{\textrm{lex}}i\}$ and for $k\geq 2$, 
$V_{i}^{k}=V_{i}^{1}\cap\{j\in\Z^{d}\,;\,\vert i-j\vert\geq k\}$. For all $n$ in $\N\backslash\{0\}$ and all $k$ in $\{1,...,\vert\Lambda_n\vert\}$, we define
$$
\textrm{E}_{k}^{(n)}=\varphi(\{1,..,k\})\cap V_{\varphi(k)}^{M_n}\quad\textrm{and}\quad
S_{\varphi(k)}^{M_n}(Y)=\sum_{i\in\textrm{E}_{k}^{(n)}}Y_{i}.
$$
For all function $h$ from $\R$ to $\R$, we define $h_{k-1,l}^{M_n}=h\left(S_{\varphi(k)}^{M_n}(Y)+S_{\varphi(l)}^c(\gamma)\right)$. 
Our aim is to show that
\begin{equation}\label{equation1bis}
\lim_{n\to\infty}\sum_{k=1}^{\vert\Lambda_n\vert}\E\left(Y_{\varphi(k)}h_{k-1,k+1}^{'}-Y_{\varphi(k)}\left(S_{\varphi(k-1)}(Y)-S_{\varphi(k)}^{M_n}(Y)\right)h_{k-1,k+1}^{''}\right)=0.
\end{equation}
First, we use the decomposition
$$
Y_{\varphi(k)}h_{k-1,k+1}^{'}=Y_{\varphi(k)}h_{k-1,k+1}^{'M_n}+Y_{\varphi(k)}\left(h_{k-1,k+1}^{'}-h_{k-1,k+1}^{'M_n}\right).
$$
Applying again Taylor's formula,
$$
Y_{\varphi(k)}(h_{k-1,k+1}^{'}-h_{k-1,k+1}^{'M_n})=Y_{\varphi(k)}\left(S_{\varphi(k-1)}(Y)-S_{\varphi(k)}^{M_n}(Y)\right)h_{k-1,k+1}^{''}+R_{k}^{'},
$$
where 
$$
\vert R_{k}^{'}\vert\leq 2\left\vert Y_{\varphi(k)}\left(S_{\varphi(k-1)}(Y)-S_{\varphi(k)}^{M_n}(Y)\right)\left(1\wedge\vert S_{\varphi(k-1)}(Y)-S_{\varphi(k)}^{M_n}(Y)\vert\right)\right\vert.
$$ 
Since $(Y_i)_{i\in\Z^d}$ is $M_n$-dependent, we have $\E \left(Y_{\varphi(k)}h_{k-1,k+1}^{'{M_n}}\right)=0$ 
and consequently $(\ref{equation1bis})$ holds if and only if $\lim_{n\to\infty}\sum_{k=1}^{\vert\Lambda_n\vert}\E \vert R_{k}^{'}\vert=0$. 
In fact, considering the sets $W_n=\{-M_n+1,...,M_n-1\}^d$ and $W_n^{\ast}=W_n\backslash\{0\}$, it follows that
\begin{align*}
\sum_{k=1}^{\vert\Lambda_n\vert}\E \vert R_{k}^{'}\vert 
&\leq 2\E \left(\vert\overline{\Delta}_{0}\vert\left(\sum_{i\in W_n^{\ast}}\vert\overline{\Delta}_{i}\vert\right)
\left(1\wedge\frac{1}{\vert\Lambda_n\vert^{1/2}}\sum_{i\in W_n^{\ast}}\vert\overline{\Delta}_{i}\vert\right)\right)\\
&\leq2M_n^d\sup_{i\in\Z^d\backslash\{0\}}
\E(\vert\overline{\Delta}_0\overline{\Delta}_i\vert)\\
&=o(1)\qquad\textrm{(by Lemma \ref{lemme-technique})}.
\end{align*}
In order to obtain ($\ref{equation1}$) it remains to control
$$
\textrm{F}_{1}=\E \left(\sum_{k=1}^{\vert\Lambda_n\vert}h_{k-1,k+1}^{''}\left(\frac{Y_{\varphi(k)}^2}{2}+Y_{\varphi(k)}\left(S_{\varphi(k-1)}(Y)-S_{\varphi(k)}^{M_n}(Y)\right)-
\frac{\eta}{2\vert\Lambda_n\vert}\right)\right).
$$
Applying again Lemma $\ref{lemme-technique}$, we have
\begin{align*}
\textrm{F}_1&\leq\left\vert\E\left(\frac{1}{\vert\Lambda_n\vert}
\sum_{k=1}^{\vert\Lambda_n\vert}h_{k-1,k+1}^{''}\left(\overline{\Delta}_{\varphi(k)}^2-\E(\overline{\Delta}_0^2)\right)\right)\right\vert
+\left\vert\eta-\E\left(\overline{\Delta}_0^2\right)\right\vert+2\sum_{j\in V_0^1\cap W_n}\E\vert\overline{\Delta}_0\overline{\Delta}_{j}\vert\\
&\leq\left\vert\E\left(\frac{1}{\vert\Lambda_n\vert}
\sum_{k=1}^{\vert\Lambda_n\vert}h_{k-1,k+1}^{''}\left(\overline{\Delta}_{\varphi(k)}^2-\E(\overline{\Delta}_0^2)\right)\right)\right\vert
+o(1).
\end{align*}
So, it suffices to prove that
$$
\textrm{F}_2=\left\vert\E\left(\frac{1}{\vert\Lambda_n\vert}
\sum_{k=1}^{\vert\Lambda_n\vert}h_{k-1,k+1}^{''}\left(\overline{\Delta}_{\varphi(k)}^2-\E(\overline{\Delta}_0^2)\right)\right)\right\vert
$$
goes to zero as $n$ goes to infinity. In fact, we have
$\textrm{F}_2\leq\frac{1}{\vert\Lambda_n\vert}\sum_{k=1}^{\vert\Lambda_n\vert}\left(\textrm{J}_k^{(1)}(n)+\textrm{J}_k^{(2)}(n)\right)$ where 
$\textrm{J}_k^{(1)}(n)=\left\vert\E\left(h_{k-1,k+1}^{''M_n}\left(\overline{\Delta}_{\varphi(k)}^2-\E\left(\overline{\Delta}_0^2\right)\right)\right)\right\vert=0$ 
since $h_{k-1,k+1}^{''M_n}$ and $\overline{\Delta}_{\varphi(k)}$ are independent. Moreover,
\begin{align*}
\textrm{J}_k^{(2)}(n)&=\left\vert\E\left(\left(h_{k-1,k+1}^{''}-h_{k-1,k+1}^{''M_n}\right)
\left(\overline{\Delta}_{\varphi(k)}^2-\E\left(\overline{\Delta}_0^2\right)\right)\right)\right\vert\\
&\leq\E\left(\left(2\wedge\sum_{\substack{\vert i\vert<M_n \\ i\neq 0}}\frac{\vert\overline{\Delta}_i\vert}{\vert\Lambda_n\vert^{1/2}}\right)\overline{\Delta}_0^2\right)\\
&\leq\frac{1}{\sqrt{\vert\Lambda_n\vert b_n}}\,\E\left(\vert\overline{\Delta}_0\vert\sqrt{b_n}
\times\sum_{\substack{\vert i\vert<M_n \\ i\neq 0}}\vert\overline{\Delta}_0\overline{\Delta}_i\vert\right)\\
&=o(1)
\end{align*}
since $(\vert\overline{\Delta}_0\vert\sqrt{b_n})_{n\geq 1}$ is uniformly integrable and $\sum_{\substack{\vert i\vert<M_n \\ i\neq 0}}\E\vert\overline{\Delta}_0\overline{\Delta}_i\vert=o(1)$ by Lemma \ref{lemme-technique}. The proof of Theorem $\ref{convergence-loi}$ is complete.
\subsection{Proof of Theorem $\textbf{\ref{Berry-Esseen_type_clt}}$}
Let $n$ be a fixed positive integer and let $x$ be fixed in $\R$. We have $U_n(x)=\overline{U}_n(x)+R_n(x)$ where 
$$
\overline{U}_n(x)=\frac{\sqrt{\vert\Lambda_n\vert b_n}\left(\overline{f}_n(x)-\E\overline{f}_n(x)\right)}{\sqrt{f(x)\int_{\R}\K^2(t)dt}}
\quad\textrm{and}\quad R_n(x)=\frac{\sqrt{\vert\Lambda_n\vert b_n}\left(f_n(x)-\overline{f}_n(x)\right)}{\sqrt{f(x)\int_{\R}\K^2(t)dt}}.
$$ 
Denote $\overline{D}_n(x)=\sup_{t\in\R}\vert\P(\overline{U}_n(x)\leq t)-\Phi(t)\vert$ and let $p\geq 2$ be fixed. 
Arguing as in Theorem 2.2 in \cite{Elmachkouri2010}, we have
\begin{equation}\label{inequality_vitesse_Un_via_Un_barre}
D_n(x)\leq \overline{D}_n(x)+\|R_n\|_{p}^{\frac{p}{p+1}}.
\end{equation}
Denoting $\sigma^2=f(x)\int_{\R}\K^2(t)dt$ and $\sigma_n^2=\E\left(\overline{U}_n^2\right)$, we have
\begin{align*}
\overline{D}_n(x)&=\sup_{t\in\R}\vert\P(\overline{U}_n(x)\leq t)-\Phi(t)\vert\\
&\leq \sup_{t\in\R}\vert\P(\overline{U}_n(x)\leq t)-\Phi\left(t/\sigma_n\right)\vert+
\sup_{t\in\R}\vert \Phi\left(t/\sigma_n\right)-\Phi\left(t\right)\vert\\
&=\sup_{t\in\R}\vert\P(\overline{U}_n(x)\leq t\sigma_n)-\Phi\left(t\right)\vert+
\sup_{t\in\R}\vert \Phi\left(t/\sigma_n\right)-\Phi\left(t\right)\vert.
\end{align*}
Applying the Berry-Esseen's type theorem for $m_n$-dependent random fields established by Chen and Shao (\cite{Chen--Shao2004}, Theorem 2.6), we obtain
\begin{equation}\label{Berry-Esseen_bound_for_mn_dependent_random_fields}
\sup_{t\in\R}\vert\P(\overline{U}_n(x)\leq t\sigma_n)-\Phi\left(t\right)\vert\leq \frac{\kappa\int_{\R}\vert K(t)\vert^\tau f(x-tb_n)dt\,\,m_n^{(\tau-1)d}}{\sigma^\tau(\vert\Lambda_n\vert b_n)^{\frac{\tau}{2}-1}}.
\end{equation}
Arguing as in Yang et al. (\cite{Yang_et_al2012}, p. 456), we have 
\begin{align*}
\sup_{t\in\R}\vert \Phi\left(t/\sigma_n\right)-\Phi\left(t\right)\vert
&\leq (2\pi e)^{-\frac{1}{2}}(\sigma_n-1)\ind{\sigma_n\geq 1}+(2\pi e)^{-\frac{1}{2}}\left(\frac{1}{\sigma_n}-1\right)\ind{0<\sigma_n< 1} \\
&\leq(2\pi e)^{-\frac{1}{2}}\max\left\{\vert\sigma_n-1\vert,\frac{\vert\sigma_n-1\vert}{\sigma_n}\}\right\}\\
&\leq\kappa \max\left\{\vert\sigma_n-1\vert,\frac{\vert\sigma_n-1\vert}{\sigma_n}\}\right\}\times(\sigma_n+1)\\
&\leq\kappa\vert\sigma_n^2-1\vert.
\end{align*}
So, we derive 
\begin{equation}\label{bound_for_Dn_barre}
\overline{D}_n(x)\leq 
\frac{\kappa\int_{\R}\vert K(t)\vert^\tau f(x-tb_n)dt\,\,m_n^{(\tau-1)d}}{\sigma^\tau(\vert\Lambda_n\vert b_n)^{\frac{\tau}{2}-1}}+\kappa\vert\sigma_n^2-1\vert.
\end{equation}
Using (\ref{developpement_norme_2_au_carre}), we have also
\begin{equation}\label{borne_sigma_2_moins_un}
\vert\sigma_n^2-1\vert
\leq\frac{1}{\sigma^2}\left\vert\E(\overline{Z}_0^2(x))-\sigma^2\right\vert
+\sum_{\substack{j\in\Z^d\backslash\{0\} \\ \vert j\vert<M_n}}
\left\vert\E\left(\overline{Z}_0(x)\overline{Z}_j(x)\right)
\right\vert.
\end{equation}
Noting that $\|\K_0(x)\|_1=O(b_n)$ and $\|\K_0(x)\|_2=O(\sqrt{b_n})$ and using the following lemma,
\begin{lemma}\label{controle_norme_de_la_difference_de_K_0_et_K_0_barre}
For all $p>1$, any positive integer $n$ and any $x$ in $\R$,
$$
\|\emph{K}_0(x)-\overline{\emph{K}}_0(x)\|_p\leq\frac{\sqrt{2p}}{b_n}\sum_{\vert j\vert>m_n}\delta_{j,p},
$$
\end{lemma}
we obtain
\begin{align*}
\left\vert\E(\overline{Z}_0^2(x))-\E(Z_0^2(x))\right\vert
&=\frac{1}{b_n}\left\vert\E(\overline{\K}_0^2(x))-\E(\K_0^2(x))\right\vert\\
&\leq\frac{1}{b_n}\|\K_0(x)\|_2\|\K_0(x)-\overline{\K}_0(x)\|_2\\
&\leq \frac{\kappa}{b_n^{3/2}}\sum_{\vert j\vert>m_n}\delta_j
\end{align*}
and
\begin{align*}
\left\vert\E(Z_0^2(x))-\sigma^2\right\vert
&=\left\vert\frac{1}{b_n}\left(\E(\K_0^2(x))-\left(\E(\K_0(x)\right)^2\right)-f(x)\int_{\R}\K^2(t)dt\right\vert\\
&\leq\left\vert\frac{1}{b_n}\E(\K_0^2(x))-f(x)\int_{\R}\K^2(t)dt\right\vert+\frac{1}{b_n}\left(\E(\K_0(x)\right)^2\\
&\leq \int_{\R}\K^2(v)\vert f(x-vb_n)-f(x)\vert dv+O(b_n)\\
&\leq\kappa\,b_n\int_{\R}\vert v\vert\K^2(v)dv+O(b_n)\\
&=O(b_n).
\end{align*}
Hence,
\begin{equation}\label{control_esperance_Z_0_barre_au_carre_moins_sigma_carre}
\left\vert\E(\overline{Z}_0^2(x))-\sigma^2\right\vert
\leq\frac{\kappa}{b_n^{3/2}}\sum_{\vert j\vert>m_n}\delta_j+O(b_n).
\end{equation}
Now, let $i\neq 0$ be fixed. We have
\begin{equation}\label{bound_for_expectation_Z0_Zi_barres}
\E\vert\overline{Z}_0(x)\overline{Z}_i(x)\vert
\leq\frac{1}{b_n}\E\vert\overline{\K}_0(x)\overline{\K}_i(x)\vert
+\frac{3}{b_n}\left(\E\vert\K_0(x)\vert\right)^2.
\end{equation}
Moreover, keeping in mind that $\vert\vert \alpha\vert-\vert \beta\vert\vert\leq\vert \alpha-\beta\vert$ for all $(\alpha,\beta)$ in $\R^2$ and applying the Cauchy-Schwarz inequality, we obtain
$$
\big\vert\E\vert\overline{\K}_0(x)\overline{\K}_i(x)\vert-\E\vert \K_0(x)\K_i(x)\vert\big\vert \leq 2\|\K_0(x)\|_2\|\K_0(x)-\overline{\K}_0(x)\|_2\\
$$
and applying Lemma \ref{controle_norme_de_la_difference_de_K_0_et_K_0_barre}, we derive
\begin{equation}\label{ecart_K0Ki_K0Kibarre}
\big\vert\E\vert\overline{\K}_0(x)\overline{\K}_i(x)\vert-\E\vert \K_0(x)\K_i(x)\vert\big\vert
\leq\frac{\kappa}{\sqrt{b_n}}\sum_{\vert j\vert>m_n}\delta_j.
\end{equation}
Combining (\ref{bound_for_expectation_Z0_Zi_barres}) and (\ref{ecart_K0Ki_K0Kibarre}), we have
\begin{equation}\label{bound_for_expectation_Z0_Zi_barres_bis}
\E\vert\overline{Z}_0(x)\overline{Z}_i(x)\vert
\leq \frac{\kappa}{b_n^{3/2}}\sum_{\vert j\vert>m_n}\delta_j+
\frac{1}{b_n}\E\vert\K_0(x)\K_i(x)\vert
+\frac{3}{b_n}\left(\E\vert\K_0(x)\vert\right)^2.
\end{equation}
Using Assumption $\textbf{(B3)}$, we obtain
\begin{align*}
\E\big\vert \textrm{K}_0(x)\textrm{K}_i(x)\big\vert
&=\iint_{\R^2}\left\vert \textrm{K}\left(\frac{x-u}{b_n}\right)\textrm{K}\left(\frac{x-v}{b_n}\right)\right\vert f_{0,i}(u,v)dudv\\
&\leq \kappa b_n^2\left(\int_{\R}\vert \textrm{K}(w)\vert dw\right)^2.
\end{align*}
Since $\E\vert\K_0(x)\vert=O(b_n)$, we derive from (\ref{bound_for_expectation_Z0_Zi_barres_bis}) that
\begin{equation}\label{bound_sum_expectation_Z0_Zi_barres}
\sum_{\substack{j\in\Z^d\backslash\{0\} \\ \vert j\vert<M_n}}
\left\vert\E\left(\overline{Z}_0(x)\overline{Z}_j(x)\right)
\right\vert
\leq \frac{\kappa M_n^d}{b_n^{3/2}}\sum_{\vert j\vert>m_n}\delta_j+O(M_n^db_n).
\end{equation}
Finally, combining (\ref{bound_for_Dn_barre}), (\ref{borne_sigma_2_moins_un}), (\ref{control_esperance_Z_0_barre_au_carre_moins_sigma_carre}) and (\ref{bound_sum_expectation_Z0_Zi_barres}), for all $\alpha>1$, we obtain
\begin{equation}\label{bound_for_Dn_barre_bis}
\overline{D}_n(x)\leq 
\frac{\kappa m_n^{d(\tau-1)}}{\sigma^\tau(\vert\Lambda_n\vert b_n)^{\frac{\tau}{2}-1}}+\frac{\kappa }{m_n^{d(\alpha-1)}b_n^{3/2}}\sum_{\vert j\vert>m_n}\vert j\vert^{d\alpha}\delta_j+O(m_n^db_n).
\end{equation}
Since there exist $\alpha>1$ and $p\geq 2$ such that 
$\sum_{i\in\Z^d}\vert i\vert^{d\alpha}\delta_{i,p}<\infty$, we derive from Lemma \ref{moment_inequality} that
\begin{equation}\label{controle_R_n_dans_Lp}
\|R_n(x)\|_p\leq\frac{\kappa \sqrt{p}}{\sigma m_n^{d(\alpha-1)}b_n^{3/2}}\sum_{i\in\Z^d}\vert i\vert^{d\alpha}\delta_{i,p}.
\end{equation}
Combining (\ref{inequality_vitesse_Un_via_Un_barre}), (\ref{bound_for_Dn_barre_bis}) and (\ref{controle_R_n_dans_Lp}), we obtain
\begin{equation}\label{majoration_finale_Dnx}
D_n(x)\leq\kappa\left(m_n^{d(\tau-1)}\left(b_n+\frac{1}{(\vert\Lambda_n\vert b_n)^{\frac{\tau}{2}-1}}\right)+\left(\frac{1}{m_n^{d(\alpha-1)}b_n^{3/2}}\right)^{\frac{p}{p+1}}\right)
\end{equation}
for all $2<\tau\leq 3$, all $p\geq 2$ and all $\alpha>1$ such that $\sum_{i\in\Z^d}\vert i\vert^{d\alpha}\delta_{i,p}<\infty$. Optimizing in $m_n$ we derive
$$
D_n(x)\leq\kappa\,b_n^{\theta_1}\left(b_n+\frac{1}{(\vert\Lambda_n\vert b_n)^{\frac{\tau}{2}-1}}\right)^{\theta_2}
$$
where 
$$
\theta_1=\frac{3p(1-\tau)}{2(\tau-1)(p+1)+2p(\alpha-1)}
\quad\textrm{and}\quad
\theta_2=\frac{p(\alpha-1)}{(\tau-1)(p+1)+p(\alpha-1)}.
$$
Finally, choosing $b_n=\vert\Lambda_n\vert^{\frac{2}{\tau}-1}$, we obtain
$D_n(x)\leq\kappa\vert\Lambda_n\vert^{-\theta}$ where 
$$
\theta=\left(\frac{1}{2}-\frac{1}{\tau}\right)\frac{3p(1-\tau)+2p(\alpha-1)}{(\tau-1)(p+1)+p(\alpha-1)}.
$$
The proof of Theorem \ref{Berry-Esseen_type_clt} is complete.
\section{Appendix}
{\em Proof of Lemma $\ref{mn}$}. We follow the proof by Bosq et al. (\cite{Bosq--Merlevede--Peligrad1999}, pages 88-89). 
First, $m_n$ goes to infinity since $v_n=\big[b_n^{-\frac{1}{2d}}\big]$ goes to infinity and $m_n\geq v_n$. For all positive integer $m$, we consider $r(m)=\sum_{\vert i\vert>m}\vert i\vert^{\frac{5d}{2}}\,\delta_i$. Since $\textbf{(A4)}$ holds, $r(m)$ converges to zero as $m$ goes to infinity. Moreover, $m_n^db_n\leq\max\left\{\sqrt{b_n},\kappa\left(r(v_n)^{1/3}+b_n\right)\right\}\converge{n}{\infty}{ }0$ and $m_n^d\geq\frac{1}{b_n}\left(r\left(v_n\right)\right)^{1/3}\geq\frac{1}{b_n}\left(r\left(m_n\right)\right)^{1/3}$ since $v_n\leq m_n$. Finally, we obtain
$$
\frac{1}{\left(m_n^db_n\right)^{3/2}}\sum_{\vert i\vert>m_n}\vert i\vert^{\frac{5d}{2}}\,\delta_i\leq\sqrt{r(m_n)}\converge{n}{\infty}{ }0.
$$
The proof of Lemma $\ref{mn}$ is complete.\\
\\
{\em Proof of Lemma $\ref{moment_inequality}$}. Let $p>1$ be fixed. We follow the proof of Proposition 1 in \cite{Elmachkouri--Volny--Wu2013}. For all $i$ in $\Z^d$ and all $x$ in $\R$, we denote 
$R_i=\K_i(x)-\overline{\K}_i(x)$. Since there exists a measurable function $\textrm{H}$ such that 
$R_i=\textrm{H}(\varepsilon_{i-s};s\in\Z^d)$, we are able to define the physical dependence measure coefficients $(\delta^{(n)}_{i,p})_{i\in\Z^d}$ 
associated to the random field $(R_i)_{i\in\Z^d}$. We recall that $\delta^{(n)}_{i,p}=\|R_i-R_i^{\ast}\|_p$ where 
$R_i^{\ast}=\textrm{H}(\varepsilon^{\ast}_{i-s};s\in\Z^d)$ and $\varepsilon_j^{\ast}=\varepsilon_j\ind{\{j\neq 0\}}+\varepsilon_0^{'}\ind{\{j=0\}}$ 
for all $j$ in $\Z^d$. In other words, we obtain $R_i^{\ast}$ from $R_i$ by just replacing $\varepsilon_0$ by its copy $\varepsilon_0^{'}$ (see \cite{Wu2005}). 
Let $\tau:\Z\to\Z^d$ be a bijection. For all $l\in\Z$, for all $i\in\Z^d$, we denote $P_lR_i:=\E(R_i\vert\F_l)-\E(R_i\vert\F_{l-1})$ where $\F_l=\sigma\left(\varepsilon_{\tau(s)};s\leq l\right)$ and $R_i=\sum_{l\in\Z}P_lR_i$. Consequently, 
$\left\|\sum_{i\in\Lambda_n}a_iR_i\right\|_p=\left\|\sum_{l\in\Z}\sum_{i\in \Lambda_n}a_iP_lR_i\right\|_p$ and
applying the Burkholder inequality (cf. \cite{Hall--Heyde1980}, page 23) for the martingale difference sequence 
$\left(\sum_{i\in \Lambda_n}a_iP_lR_i\right)_{l\in\Z}$, we obtain
\begin{equation}\label{Burkholder_inequality}
\left\|\sum_{i\in\Lambda_n}a_iR_i\right\|_p
\leq\left(2p\sum_{l\in\Z}\left\|\sum_{i\in \Lambda_n}a_iP_lR_i\right\|_p^2\right)^{\frac{1}{2}}
\leq\left(2p\sum_{l\in\Z}\left(\sum_{i\in \Lambda_n}\vert a_i\vert\left\|P_lR_i\right\|_p\right)^2\right)^{\frac{1}{2}}.
\end{equation}
Moreover, by the Cauchy-Schwarz inequality, we have
\begin{equation}\label{Cauchy-Schwarz_inequality}
\left(\sum_{i\in \Lambda_n}\vert a_i\vert\left\|P_lR_i\right\|_p\right)^2
\leq\sum_{i\in \Lambda_n}a_i^2\left\|P_lR_i\right\|_p\times\sum_{i\in\Lambda_n}\|P_lR_i\|_p.
\end{equation}
Let $l$ in $\Z$ and $i$ in $\Z^d$ be fixed.
$$
\left\|P_lR_i\right\|_p=\left\|\E(R_i\vert\F_l)-\E(R_i\vert\F_{l-1})\right\|_p=\left\|\E(R_0\vert T^i\F_l)-\E(R_0\vert T^i\F_{l-1})\right\|_p
$$
where $T^i\F_l=\sigma\left(\varepsilon_{\tau(s)-i};s\leq l\right)$. Hence,
\begin{align*}
\left\|P_lR_i\right\|_p&=\left\|\E\left(\textrm{H}\left((\varepsilon_{-s})_{s\in\Z^d}\right)\vert T^i\F_l\right)-
\E\left(\textrm{H}\left((\varepsilon_{-s})_{s\in\Z^d\backslash\{i-\tau(l)\}};\varepsilon^{'}_{\tau(l)-i}\right)\vert T^i\F_{l}\right)\right\|_p\\
&\leq \left\|\textrm{H}\left((\varepsilon_{-s})_{s\in\Z^d}\right)-\textrm{H}\left((\varepsilon_{-s})_{s\in\Z^d\backslash\{i-\tau(l)\}};\varepsilon^{'}_{\tau(l)-i}\right)\right\|_p\\
&=\left\|\textrm{H}\left((\varepsilon_{i-\tau(l)-s})_{s\in\Z^d}\right)-\textrm{H}\left((\varepsilon_{i-\tau(l)-s})_{s\in\Z^d\backslash\{i-\tau(l)\}};\varepsilon^{'}_{0}\right)\right\|_p\\
&=\left\|R_{i-\tau(l)}-R_{i-\tau(l)}^{\ast}\right\|_p\\
&=\delta^{(n)}_{i-\tau(l),p}.
\end{align*}
Consequently, $\sum_{i\in\Z^d}\|P_lR_i\|_p\leq\sum_{j\in\Z^d}\delta^{(n)}_{j,p}$ and combining (\ref{Burkholder_inequality}) and (\ref{Cauchy-Schwarz_inequality}), we obtain
$$
\left\|\sum_{i\in\Lambda_n}a_iR_i\right\|_p\leq\left(2p\sum_{j\in\Z^d}\delta^{(n)}_{j,p}\sum_{i\in \Lambda_n}a_i^2\sum_{l\in\Z}\left\|P_lR_i\right\|_p\right)^{\frac{1}{2}}.
$$
Similarly, for all $i$ in $\Z^d$, we have $\sum_{l\in\Z}\|P_lR_i\|_p\leq\sum_{j\in\Z^d}\delta^{(n)}_{j,p}$ and we derive
\begin{equation}\label{moment_inequality_temp}
\left\|\sum_{i\in\Lambda_n}a_iR_i\right\|_p\leq\left(2p\sum_{i\in \Lambda_n}a_i^2\right)^{\frac{1}{2}}\sum_{i\in\Z^d}\delta^{(n)}_{i,p}.
\end{equation}
Since $\overline{\K}_i^{\ast}=\E\left(\K_i^{\ast}(x)\big\vert\F_{n,i}^{\ast}\right)$ where 
$\F_{n,i}^{\ast}=\sigma\left(\varepsilon^{\ast}_{i-s}\,;\,\vert s\vert\leq m_n\right)$ and 
$\left(\K_i(x)-\overline{\K}_i(x)\right)^{\ast}=\K_i^{\ast}(x)-\overline{\K}_i^{\ast}(x)$, 
we derive $\delta_{i,p}^{(n)}\leq 2\|\K_i(x)-\K_i^{\ast}(x)\|_p$. Since $\K$ is Lipschitz, we obtain 
\begin{equation}\label{delta_i_p_n_inequality1}
\delta_{i,p}^{(n)}\leq\frac{2\delta_{i,p}}{b_n}
\end{equation}
where $\delta_{i,p}=\|X_i-X^{\ast}_i\|_p$. Morever, we have also 
\begin{equation}\label{delta_i_p_n_inequality1bis}
\delta_{i,p}^{(n)}\leq 2\|\K_0(x)-\overline{\K}_0(x)\|_p.
\end{equation}
Combining (\ref{delta_i_p_n_inequality1bis}) and Lemma \ref{controle_norme_de_la_difference_de_K_0_et_K_0_barre}, we derive
\begin{equation}\label{delta_i_p_n_inequality2}
\delta_{i,p}^{(n)}\leq\frac{\sqrt{8p}}{b_n}\sum_{\vert j\vert>m_n}\delta_{j,p}.
\end{equation}
Combining ($\ref{delta_i_p_n_inequality1}$) and ($\ref{delta_i_p_n_inequality2}$), we obtain 
$$
\sum_{i\in\Z^d}\delta_{i,p}^{(n)}\leq \frac{m_n^d\sqrt{8p}}{b_n}\sum_{\vert j\vert>m_n}\delta_{j,p}+\frac{2}{b_n}\sum_{\vert j\vert>m_n}\delta_{j,p}
\leq\frac{2\sqrt{8p}m_n^d}{b_n}\sum_{\vert j\vert>m_n}\delta_{j,p}.
$$
The proof of Lemma \ref{moment_inequality} is complete.\\
\\
{\em Proof of Lemma $\ref{lemme_technique_0}$}.  Let $s$ and $t$ be fixed in $\R$. Since $\E\left(\overline{\K}_0(s)\overline{\K}_0(t)\right)=\E\left(\K_0(s)\overline{\K}_0(t)\right)$, we have 
$$
\big\vert\E\left(\overline{\K}_0(s)\overline{\K}_0(t)\right)-\E\left(\K_0(s)\K_0(t)\right)\big\vert
\leq \|\K_0(s)\|_2\|\K_0(t)-\overline{\K}_0(t)\|_2.
$$
Keeping in mind that $\|\K_0(s)\|_2=O(\sqrt{b_n})$ and using Lemma \ref{controle_norme_de_la_difference_de_K_0_et_K_0_barre}, we have
$$
\big\vert\E\left(\overline{\K}_0(s)\overline{\K}_0(t)\right)-\E\left(\K_0(s)\K_0(t)\right)\big\vert
\leq \frac{\kappa}{\sqrt{b_n}}\sum_{\vert j\vert>m_n}\delta_j.
$$
Since $b_n\vert\E(Z_0(s)Z_0(t))-\E(\overline{Z}_0(s)\overline{Z}_0(t)\vert=\vert\E\left(\K_0(s)\K_0(t)\right)-\E\left(\overline{\K}_0(s)\overline{\K}_0(t)\right)\vert$, we have
\begin{equation}\label{ecart_Z_0_Z_0_barre}
M_n^d\vert\E(Z_0(s)Z_0(t))-\E(\overline{Z}_0(s)\overline{Z}_0(t)\vert
\leq\frac{\kappa}{(m_n^db_n)^{3/2}}\sum_{\vert j\vert>m_n}\vert j\vert^{\frac{5d}{2}}\,\delta_j.
\end{equation}
Moreover, keeping in mind Assumptions $\textbf{(A1)}$, $\textbf{(A2)}$ and $\textbf{(A4)}$, we have
\begin{equation}\label{limit_f_sigma2}
\lim_{n}\frac{1}{b_n}\E\left(\K_0(s)\K_0(t)\right)
=\lim_{n}\int_{\R}\K \left(v\right)\K \left(v+\frac{t-s}{b_n}\right)f(s-vb_n)dv
=u(s,t)\,f(s)\,\int_{\R}\K^2(u)du
\end{equation}
where $u(s,t)=1$ if $s=t$ and $u(s,t)=0$ if $s\neq t$. We have also
\begin{equation}\label{limit_zero}
\lim_{n}\frac{1}{b_n}\E\K_0(s)\E\K_0(t)=\lim_{n}b_n\int_{\R}\K (v)f(s-vb_n)dv\int_{\R}\K (w)f(t-wb_n)dw=0.
\end{equation}
Let $x$ be fixed in $\R$. Choosing $s=t=x$ and combining (\ref{ecart_Z_0_Z_0_barre}), (\ref{limit_f_sigma2}), 
(\ref{limit_zero}) and Lemma \ref{mn}, we obtain 
$\E(\overline{Z}_0^2(x))$ goes to $f(x)\int_{\R}\K^2(u)du$ as $n$ goes to infinity.\\
In the other part, let $i\neq 0$ be fixed in $\Z^d$ and let $s$ and $t$ be fixed in $\R$. We have
\begin{equation}\label{borne_Z_0_Z_i}
\E\vert\overline{Z}_0(s)\overline{Z}_i(t)\vert\leq\frac{1}{b_n}\E\big\vert \overline{\K}_0(s)\overline{\K}_i(t)\big\vert+\frac{3}{b_n}\E\big\vert\K_0(s)\big\vert\,\E\big\vert \K_0(t)\big\vert.
\end{equation}
Keeping in mind that $\vert\vert \alpha\vert-\vert \beta\vert\vert\leq\vert \alpha-\beta\vert$ for all $(\alpha,\beta)$ in $\R^2$ and applying the Cauchy-Schwarz inequality, we obtain 
\begin{equation}\label{by_Schwarz_inequality}
\big\vert\E\vert\overline{\K}_0(s)\overline{\K}_i(t)\vert-\E\vert \K_0(s)\K_i(t)\vert\big\vert
\leq \|\overline{\K}_0(s)\|_2\|\overline{\K}_0(t)-\K_0(t)\|_2+\|\K_0(t)\|_2\|\overline{\K}_0(s)-\K_0(s)\|_2
\end{equation}
Applying again Lemma \ref{controle_norme_de_la_difference_de_K_0_et_K_0_barre}, we obtain
\begin{equation}\label{ecart_overline_K_O_K_i_et_K_0_K_i}
\frac{M_n^d}{b_n}\big\vert\E\vert\overline{\K}_0(s)\overline{\K}_i(t)\vert-\E\vert \K_0(s)\K_i(t)\vert\big\vert
\leq\frac{\kappa}{(m_n^db_n)^{3/2}}\sum_{\vert j\vert>m_n}\vert j\vert^{\frac{5d}{2}}\delta_j.
\end{equation}
Since Assumptions $\textbf{(A1)}$ and $\textbf{(A4)}$ hold and $M_n^db_n=o(1)$, we have
\begin{equation}\label{O_m_n_b_n}
\frac{M_n^d}{b_n}\E\big\vert\K_0(s)\big\vert\,\E\big\vert \K_0(t)\big\vert=M_n^db_n\int_{\R}\vert \K (u)\vert f(s-ub_n)du\int_{\R}\vert \K (v)\vert f(t-vb_n)dv=o(1).
\end{equation}
Moreover, using Assumption $\textbf{(B3)}$, we have
\begin{align*}
\E\big\vert \textrm{K}_0(s)\textrm{K}_i(t)\big\vert
&=\iint_{\R^2}\left\vert \textrm{K}\left(\frac{s-u}{b_n}\right)\textrm{K}\left(\frac{t-v}{b_n}\right)\right\vert f_{0,i}(u,v)dudv\\
&\leq \kappa b_n^2\left(\int_{\R}\vert \textrm{K}(w)\vert dw\right)^2.
\end{align*}
So, using again Assumption $\textbf{(A4)}$ and $M_n^db_n=o(1)$, we derive
\begin{equation}\label{O_m_n_b_n_bis}
\frac{M_n^d}{b_n}\E\big\vert \K_0(s)\K_i(t)\big\vert=o(1).
\end{equation}
Combining (\ref{borne_Z_0_Z_i}), (\ref{ecart_overline_K_O_K_i_et_K_0_K_i}), (\ref{O_m_n_b_n}), (\ref{O_m_n_b_n_bis}) and Lemma \ref{mn}, we obtain 
\begin{equation}\label{overline_Z_0_overline_Z_i}
M_n^d\sup_{i\in\Z^d\backslash\{0\}}\E\vert\overline{Z}_0(s)\overline{Z}_i(t)\vert=o(1).
\end{equation}
The proof of Lemma $\ref{lemme_technique_0}$ is complete.\\
\\
{\em Proof of Lemma $\ref{lemme-technique}$}. Let $x$ and $y$ be two distinct real numbers. Noting that 
\begin{align*}
\E(\Delta_0^2)&=\lambda_1^2\E (Z_0^2(x))+\lambda_2^2\E (Z_0^2(y))+2\lambda_1\lambda_2\E(Z_0(x)Z_0(y))\\
\E(\overline{\Delta}_0^2)&=\lambda_1^2\E (\overline{Z}_0^2(x))+\lambda_2^2\E (\overline{Z}_0^2(y))
+2\lambda_1\lambda_2\E(\overline{Z}_0(x)\overline{Z}_0(y))
\end{align*}
and using (\ref{ecart_Z_0_Z_0_barre}) and Lemma \ref{mn}, we obtain
\begin{equation}\label{ecart_Delta_et_overline_Delta}
 \lim_{n\to\infty}M_n^d\vert\E(\Delta_0^2)-\E(\overline{\Delta}_0^2)\vert=0.
\end{equation}
Combining (\ref{limit_f_sigma2}) and (\ref{ecart_Delta_et_overline_Delta}), we derive that $\E(\overline{\Delta}_0^2)$ converges to 
$\eta=\left(\lambda_1^2f(x)+\lambda_2^2f(y)\right)\int_{\R}\K^2(u)du$.\\
Let $i\neq 0$ be fixed in $\Z^d$. Combining (\ref{overline_Z_0_overline_Z_i}) and
\begin{equation}\label{delta_O_delta_i}
\E\vert\overline{\Delta}_0\overline{\Delta}_i\vert\leq\lambda_1^2\E\vert\overline{Z}_0(x)\overline{Z}_i(x)\vert+\lambda_2^2\E\vert\overline{Z}_0(y)\overline{Z}_i(y)\vert
+\lambda_1\lambda_2\E\vert\overline{Z}_0(x)\overline{Z}_i(y)\vert+\lambda_1\lambda_2\E\vert\overline{Z}_0(y)\overline{Z}_i(x)\vert,
\end{equation}
we obtain $M_n^d\sup_{i\in\Z^d\backslash\{0\}}\E\vert\overline{\Delta}_0\overline{\Delta}_i\vert=o(1)$. The proof of Lemma \ref{lemme-technique} is complete.\\
\\
{\em Proof of Lemma $\ref{controle_norme_de_la_difference_de_K_0_et_K_0_barre}$}. Let $p>1$ be fixed. 
We consider the sequence $(\Gamma_n)_{n\geq 0}$ of finite subsets of $\Z^d$ defined by 
$\Gamma_0=\{(0,...,0)\}$ and for all $n$ in $\N\backslash\{0\}$, $\Gamma_n=\{i\in\Z^d\,;\,\vert i\vert=n\}$. 
For all integer $n$, let $a_n=\sum_{j=0}^n\vert\Gamma_j\vert$ and let $\tau:\N\backslash\{0\}\to\Z^d$ be the bijection defined by $\tau(1)=(0,...,0)$ and 
\begin{itemize}
\item for all $n$ in $\N\backslash\{0\}$, if $l\in\left]a_{n-1},a_n\right]$ then $\tau(l)\in\Gamma_n$,
\item for all $n$ in $\N\backslash\{0\}$, if $(i,j)\in\left]a_{n-1},a_n\right]^2$ and $i<j$ then $\tau(i)<_{\textrm{lex}}\tau(j)$
\end{itemize}
Let $(m_n)_{n\geq 1}$ be the sequence of positive integers defined by (\ref{def_mn}). 
For all $n$ in $\N\backslash\{0\}$, we recall that $\F_{n,0}=\sigma\left(\varepsilon_{-s}\,;\,\vert s\vert\leq m_n\right)$ (see (\ref{definition_K_i_z_overline_K_i_z})) and we consider also the 
$\sigma$-algebra $\G_{n}:=\sigma\left(\varepsilon_{\tau(j)}\,;\,1\leq j\leq n\right)$. By the definition of the bijection $\tau$, we have 
$1\leq j\leq a_n$ if and only if $\vert\tau(j)\vert\leq n$. Consequently $\G_{a_{m_n}}=\F_{n,0}$ and $\K_0(x)-\overline{\K}_0(x)=\sum_{l> a_{m_n}}D_l$ with $D_l=\E\left(\K_0(x)\vert\G_l\right)-\E\left(\K_0(x)\vert\G_{l-1}\right)$ for all $l$ in $\Z$. 
Let $p>1$ be fixed. Since $\left(D_l\right)_{l\in\Z}$ is a martingale-difference sequence, applying Burkholder's inequality (cf. \cite{Hall--Heyde1980}, page 23), we derive
$$
\|\K_0(x)-\overline{\K}_0(x)\|_p\leq\left(2p\sum_{l>a_{m_n}}\|D_l\|_p^2\right)^{1/2}.
$$
Denoting $\K^{'}_0(x)=\K\left(b_n^{-1}\left(x-g\left((\varepsilon_{-s})_{s\in\Z^d\backslash\{-\tau(l)\}};\varepsilon^{'}_{\tau(l)}\right)\right)\right)$, we obtain
\begin{align*}
\|D_l\|_p&=\|\E\left(\K_0(x)\vert \G_l\right)-\E\left(\K^{'}_0(x)\vert\G_l\right)\|_p\leq \|\K_0(x)-\K^{'}_0(x)\|_p\\
&\leq\frac{1}{b_n}\left\|g\left((\varepsilon_{-s})_{s\in\Z^d}\right)-g\left((\varepsilon_{-s})_{s\in\Z^d\backslash\{-\tau(l)\}};\varepsilon^{'}_{\tau(l)}\right)\right\|_p\\
&=\frac{1}{b_n}\left\|g\left((\varepsilon_{-\tau(l)-s})_{s\in\Z^d}\right)-g\left((\varepsilon_{-\tau(l)-s})_{s\in\Z^d\backslash\{-\tau(l)\}};\varepsilon^{'}_{0}\right)\right\|_p\\
&=\frac{1}{b_n}\left\|X_{-\tau(l)}-X_{-\tau(l)}^{\ast}\right\|_p=\frac{\delta_{-\tau(l),p}}{b_n}
\end{align*}
and finally
$$
\|\K_0(x)-\overline{\K}_0(x)\|_p\leq \frac{1}{b_n}\left(2p\sum_{l>a_{m_n}}\delta_{-\tau(l),p}^2\right)^{1/2}\leq\frac{\sqrt{2p}}{b_n}\sum_{\vert j\vert>m_n}\delta_{j,p}.
$$
The proof of Lemma \ref{controle_norme_de_la_difference_de_K_0_et_K_0_barre} is complete.\\
\\
\textbf{Acknowledgments}. The author is grateful to two anonymous referees for their careful reading and constructive comments. He is also indebted to Yizao Wang for pointing an error in the proof of a first version of Theorem \ref{Berry-Esseen_type_clt}.
\bibliographystyle{plain}
\bibliography{xbib}

\begin{thebibliography}{10}

\bibitem{Bosq--Merlevede--Peligrad1999}
D.~Bosq, F.~Merlev{\`e}de, and M.~Peligrad.
\newblock Asymptotic normality for density kernel estimators in discrete and
  continuous time.
\newblock {\em J. Multivariate Anal.}, 68(1):78--95, 1999.

\bibitem{Carbon-Hallin-Tran}
M.~Carbon, M.~Hallin, and L.T. Tran.
\newblock Kernel density estimation for random fields: the $l_{1}$ theory.
\newblock {\em Journal of nonparametric {S}tatistics}, 6:157--170, 1996.

\bibitem{Carbon-Tran-Wu}
M.~Carbon, L.T. Tran, and B.~Wu.
\newblock Kernel density estimation for random fields.
\newblock {\em Statist. Probab. Lett.}, 36:115--125, 1997.

\bibitem{Casti1985}
J.~L. Casti.
\newblock {\em Nonlinear system theory}, volume 175 of {\em Mathematics in
  Science and Engineering}.
\newblock Academic Press Inc., Orlando, FL, 1985.

\bibitem{Chen--Shao2004}
Q.~M. Chen, L. H. Y.~Shao.
\newblock Normal approximation under local dependence.
\newblock {\em Ann. of {P}robab.}, 32:1985--2028, 2004.

\bibitem{Cheng-Ho2006}
T.-L. Cheng and H.-C. Ho.
\newblock Central limit theorems for instantaneous filters of linear random
  fields on {$\bold Z^2$}.
\newblock In {\em Random walk, sequential analysis and related topics}, pages
  71--84. World Sci. Publ., Hackensack, NJ, 2006.

\bibitem{Cheng-Ho-Lu2008}
T-L. Cheng, H-C. Ho, and X.~Lu.
\newblock A note on asymptotic normality of kernel estimation for linear random
  fields on {${\bf Z}^2$}.
\newblock {\em J. Theoret. Probab.}, 21(2):267--286, 2008.

\bibitem{Dedecker1998}
J.~Dedecker.
\newblock A central limit theorem for stationary random fields.
\newblock {\em Probab. Theory Relat. Fields}, 110:397--426, 1998.

\bibitem{Dedecker_Merlevede2002}
J.~Dedecker and F.~Merlev\`ede.
\newblock Necessary and sufficient conditions for the conditional central limit
  theorem.
\newblock {\em {A}nnals of {P}robability}, 30(3):1044--1081, 2002.

\bibitem{Elmachkouri2010}
M.~El~Machkouri.
\newblock Berry-{E}sseen's central limit theorem for non-causal linear
  processes in {H}ilbert spaces.
\newblock {\em African Diaspora Journal of Mathematics}, 10(2):1--6, 2010.

\bibitem{Elmachkouri2011}
M.~El~Machkouri.
\newblock Asymptotic normality for the parzen-rosenblatt density estimator for
  strongly mixing random fields.
\newblock {\em Statistical Inference for Stochastic Processes}, 14(1):73--84,
  2011.

\bibitem{Elmachkouri--Volny--Wu2013}
M.~El~Machkouri, D.~Voln{\'y}, and W.~B. Wu.
\newblock A central limit theorem for stationary random fields.
\newblock {\em Stochastic Process. Appl.}, 123(1):1--14, 2013.

\bibitem{Hall--Heyde1980}
P.~Hall and C.~C. Heyde.
\newblock {\em Martingale limit theory and its application}.
\newblock Academic {P}ress, New {Y}ork, 1980.

\bibitem{Hallin--Lu--Tran2001}
M.~Hallin, Z.~Lu, and L.T. Tran.
\newblock Density estimation for spatial linear processes.
\newblock {\em Bernoulli}, 7:657--668, 2001.

\bibitem{Hallin-Lu-Tran-2004a}
M.~Hallin, Z.~Lu, and L.T. Tran.
\newblock Density estimation for spatial processes: the $l\sb 1$ theory.
\newblock {\em J. Multivariate Anal.}, 88(1):61--75, 2004.

\bibitem{Kulkarni1992}
P.~M. Kulkarni.
\newblock Estimation of parameters of a two-dimensional spatial autoregressive
  model with regression.
\newblock {\em Statist. Probab. Lett.}, 15(2):157--162, 1992.

\bibitem{Lindeberg}
J.~W. Lindeberg.
\newblock Eine neue {H}erleitung des {E}xponentialgezetzes in der
  {W}ahrscheinlichkeitsrechnung.
\newblock {\em Mathematische {Z}eitschrift}, 15:211--225, 1922.

\bibitem{Parzen1962}
E.~Parzen.
\newblock On the estimation of a probability density and the mode.
\newblock {\em Ann. Math. Statist.}, 33:1965--1976, 1962.

\bibitem{Ros}
M.~Rosenblatt.
\newblock A central limit theorem and a strong mixing condition.
\newblock {\em Proc. {N}at. {A}cad. {S}ci. USA}, 42:43--47, 1956.

\bibitem{Rugh1981}
W.~J. Rugh.
\newblock {\em Nonlinear system theory}.
\newblock Johns Hopkins Series in Information Sciences and Systems. Johns
  Hopkins University Press, Baltimore, Md., 1981.

\bibitem{Silverman1986}
B.W. Silverman.
\newblock {\em Density {E}stimation for {S}tatistics and {D}ata {A}nalysis}.
\newblock Chapman and {H}all, {L}ondon, 1986.

\bibitem{Tran}
L.T. Tran.
\newblock Kernel density estimation on random fields.
\newblock {\em J. Multivariate Anal.}, 34:37--53, 1990.

\bibitem{Wang--Woodroofe2011}
Y.~Wang and M.~Woodroofe.
\newblock On the asymptotic normality of kernel density estimators for causal
  linear random fields.
\newblock {\em J. Multivariate Anal.}, 123:201--213, 2014.

\bibitem{Wu2005}
W.~B. Wu.
\newblock Nonlinear system theory: another look at dependence.
\newblock {\em Proc. Natl. Acad. Sci. USA}, 102(40):14150--14154 (electronic),
  2005.

\bibitem{Wu2011}
W.~B. Wu.
\newblock Asymptotic theory for stationary processes.
\newblock {\em Statistics and Its Interface}, 0:1--20, 2011.

\bibitem{Wu-Mielniczuk}
W.B. Wu and J.~Mielniczuk.
\newblock Kernel density estimation for linear processes.
\newblock {\em Ann. Statist.}, 30:1441--1459, 2002.

\bibitem{Yang_et_al2012}
W.~Yang, X.~Wang, X.~Li, and S.~Hu.
\newblock Berry-{E}ss\'een bound of sample quantiles for {$\phi$}-mixing random
  variables.
\newblock {\em J. Math. Anal. Appl.}, 388(1):451--462, 2012.

\end{thebibliography}
\end{document}